\documentclass[11pt]{article}
\usepackage{mathrsfs}
\usepackage{amsfonts}
\usepackage{natbib}
\usepackage{amsmath}
\usepackage{amsthm}
\usepackage{amssymb}
\usepackage{color}
\usepackage{lscape}
\usepackage{rotating}
\usepackage{caption2}
\usepackage{subfigure}
\usepackage{graphicx}

\usepackage{amssymb}

\def\be{\begin{equation}}
\def\ee{\end{equation}}
\def\bea{\begin{eqnarray}}
\def\eea{\end{eqnarray}}
\def\bd{\begin{displaymath}}
\def\ed{\end{displaymath}}
\def\bda{\begin{eqnarray*}}
\def\eda{\end{eqnarray*}}
\def\bsm{\begin{small}}
\def\esm{\end{small}}

\def\t0{\theta_0}

\def\ha1{\hat \beta_1}

\newcommand{\beeta}{\mbox{\boldmath $\eta$}}
\newcommand{\bvarepsilon}{\mbox{\boldmath $\varepsilon$}}
\newcommand{\bzeta}{\mbox{\boldmath $\zeta$}}

\newcommand{\bSigma}{\mbox{\boldmath $\Sigma$}}

\def\bnt{\begin{enumerate}}
\def\ent{\end{enumerate}}
\def\T{{ \mathrm{\scriptscriptstyle T} }}

\def\bsc{\begin{scriptsize}}
\def\esc{\end{scriptsize}}

\newtheorem{tm}{Theorem}[section]

\newtheorem{la}{Lemma}[section]

\newtheorem{pn}{Proposition}[section]

\theoremstyle{definition}
\newtheorem{as}{Condition}[section]

\newcommand{\bg}{\mbox{\bf g}}
\newcommand{\bu}{\mbox{\bf u}}

\newcommand{\bM}{\mbox{\bf M}}

\newcommand{\bI}{\mbox{\bf I}}
\newcommand{\bD}{\mbox{\bf D}}
\newcommand{\bA}{\mbox{\bf A}}

\newcommand{\bH}{\mbox{\bf H}}
\newcommand{\bh}{\mbox{\bf h}}
\newcommand{\bx}{\mbox{\bf x}}
\newcommand{\bfe}{\mbox{\bf e}}
\newcommand{\bw}{\mbox{\bf w}}
\newcommand{\by}{\mbox{\bf y}}
\newcommand{\ba}{\mbox{\bf a}}
\newcommand{\bb}{\mbox{\bf b}}

\newcommand{\bz}{\mbox{\bf z}}
\newcommand{\bfeta}{\mbox{\boldmath $\eta$}}
\newcommand{\bve}{\mbox{\boldmath $\varepsilon$}}

\newcommand{\calM}{\mathcal{M}}
\newcommand{\calF}{\mathcal{F}}

\newcommand{\bR}{\mbox{\bf R}}

\newcommand{\wh}{\widehat}
\newcommand{\wt}{\widetilde}

\makeatletter
\newcommand{\figcaption}{\def\@captype{figure}\caption}
\newcommand{\tabcaption}{\def\@captype{table}\caption}
\makeatother

\def\vv{1.2}
\renewcommand{\baselinestretch}{\vv}

\title{ High Dimensional Stochastic Regression with Latent Factors,
Endogeneity and Nonlinearity}

\author{Jinyuan Chang\footnote{Department of Mathematics and Statistics,
The University of Melbourne, Parkville, VIC, Australia 3010. Email:
jinyuan.chang@unimelb.edu.au. This work was completed while the first
author was a PhD student at Guanghua School of Management at Peking University.} ~~~~~~~~~~~~~~~~~~~~~~ Bin Guo\footnote{Guanghua School of Management, Peking
University, Beijing, China 100871. Email: guobin1987@pku.edu.cn.}
~~~~~~~~~~~~~~~~~~~~~~~~Qiwei Yao\footnote{Corresponding Author:
Department of Statistics, London School of Economics,
   London, WC2A 2AE, U.K. and Guanghua School of Management, Peking University, Beijing, China. Phone: +44 (0)20 7955 6767. Fax: +44 (0)20 7955 7416. Email: q.yao@lse.ac.uk. This work was partially supported by an EPSRC research
   grant.}~~~~~~~~\\[2ex]
The University of Melbourne \qquad Peking University \qquad London School of Economics }

\begin{document}

\maketitle
%\thispagestyle{empty}
%\begin{abstract}
%\begin{Large}

%\end{Large}
%\bigskip

\begin{abstract}
We consider a multivariate time series model which represents a high
dimensional vector process as a sum of three terms: a linear
regression of some observed regressors, a linear combination of some
latent and serially correlated factors, and a vector white noise. We
investigate the inference without imposing stationary conditions on
the target multivariate time series, the regressors and the
underlying factors. Furthermore we deal with the the endogeneity
that there exist correlations between the observed regressors and
the unobserved factors. We also consider the model with nonlinear
regression term which can be approximated by a linear regression
function with a large number of regressors. The convergence rates
for the estimators of regression coefficients, the number of
factors, factor loading space and factors are established under the
settings when the dimension of time series and the number of
regressors may both tend to infinity together with the sample size.
The proposed method is illustrated with both simulated and real data
examples.

\end{abstract}

\renewcommand{\baselinestretch}{\vv}

\renewcommand{\baselinestretch}{\vv}
\normalsize

\newpage
% \setcounter{page}{1}
%=====================================
\section{Introduction}

In this modern information age, the availability of large or vast
time series data bring the opportunities with challenges to time
series analysts. The demand of modelling and forecasting
high-dimensional time series arises from various practical problems
such as panel study of economic, social and natural (such as
weather) phenomena, financial market analysis, communications
engineering. On the other hand, modelling multiple time series even
with moderately large dimensions is always a challenge.  Although a
substantial proportion of the methods and the theory for univariate
autoregressive and moving average (ARMA) models has found the
multivariate counterparts, the usefulness of unregularized multiple
ARMA models suffers from the overparametrization and the lack of the
identification \citep{Lutkepohl_(2006).}. Various methods have been
developed to reduce the number of parameters and to eliminate the
non-identification issues.  For example, \cite{TiaoTsay_1989}
proposed to represent a multiple series in terms of several scalar
component models based on canonical correlation analysis,
\cite{JakemanSteeleYoung_1980} adopted a two stage regression
strategy based on instrumental variables to avoid using moving
average explicitly. Another popular approach is to represent
multiple time series in terms of a few factors defined in various
ways; see, among others, \cite{StockWatson_2005},
\cite{BaiNg_Econometrica_2002}, \cite{Fornietal_2005},
\cite{LamYaoBathia_Biometrika_2011}, and \cite{LamYao_AOS_2012}.
\cite{Davis2012} proposed a vector autoregressive (VAR) model with
sparse coefficient matrices based on partial spectral coherence.
LASSO regularization has also been applied in VAR modelling; see,
for example, \cite{ShojaieMichailidis_2010} and
\cite{SongBickel_2011}.

This paper can be viewed as a further development of
\cite{LamYaoBathia_Biometrika_2011} and \cite{LamYao_AOS_2012} which
express a high-dimensional vector time series as a linear
transformation of a low-dimensional latent factor process plus a
vector white noise. We extend their methodology and explore three
new features. We only deal with the cases when the dimension is
large in relation to the sample size. Hence all asymptotic theory is
developed when both the sample size and the dimension of time series
tend to infinity together.

Firstly, we add a regression term to the factor model. This is a useful
addition as in many applications there exist some known factors
which are among the driving forces for the dynamics of most the component
series. For example, temperature is an important factor in
forecasting household electricity consumptions. The price of a
product plays a key role in its sales over different regions. The
capital asset pricing model (CAPM) theory implies that the market index is a common factor for
pricing different assets.
% see also the example in section 6.
When the regressor and the latent factor are uncorrelated, we
estimate the regression coefficients first by the least squares
method. We then estimate the number of factors and the factor
loading space based on the residuals resulted from the regression
estimation. We show that the latter is asymptotically adaptive to
the unknown regression coefficients in the sense that the
convergence rates for estimating the factor loading space and the
factor process are the same as if the regression coefficients were
known. We also consider the models with endogeneity in the sense
that there exist correlations between the regressors and the latent
factors. We show that the factor loading space can still be identified and
estimated consistently in the presence of the endogeneity.
 However relevant instrumental variables need to be employed if the `original'
regression coefficients have to be estimated consistently. The
exploration in this direction has some overlap with
\cite{PesaranTosetti_Ecomometrica_2011}, although the models, the
inference methods and the asymptotic results in the two papers are
different.

Our second contribution lies in the fact that we do not impose
stationarity conditions on the regressors and the latent factor
process throughout the paper. This enlarges the potential
application substantially, as many important factors in practical
problems (such as temperature, calendar effects) are not stationary.
Different from the method of \cite{PanYao_Bioka_2008} which can also
handle nonstationary factors but is computationally expensive, our
approach is a direct extension of
\cite{LamYaoBathia_Biometrika_2011} and \cite{LamYao_AOS_2012} and,
hence, is applicable to the cases when the dimensions of time series
is in the order of thousands with an ordinary personal computer.

Finally, we focus on the factor models with a nonlinear regression
term. By expressing the nonlinear regression function as a linear
combination of some base functions, we turn the problem into the
model with a large number of linear regressors. Now the asymptotic
theory is established when the sample size, the dimension of time
series and the number of regressors go to infinity together.

The rest of the paper is organized as follows. Section 2 deals with
linear regression models with latent factors but without
endogeneity. The models with the endogeneity are handled in Section
3. Section 4 investigates the models with nonlinear regression term.
Simulation results are reported in Section 5. Illustration with some
stock prices included in S\&P500
% data set consisting of the 123 stock prices
is presented in Section
6. All the technical proofs are relegated to the Appendix.

\section{Regression with latent factors}

\subsection{Models}

Consider the regression model
\begin{equation}
{\bf y}_t= {\bf D}{\bf z}_t + {\bf A}{\bf x}_t +\bvarepsilon_t,
% \left\{ \begin{aligned}
%          {\bf y}_t=&~{\bf D}{\bf z}_t+\bzeta_t+\bvarepsilon_t,\\
%          \bzeta_t=&~{\bf A}{\bf x}_t,
%         \end{aligned} \right.
\label{eq:nonlinear1}
\end{equation}
where $\by_t$ and ${\bf z}_t$ are, respectively, observable $p\times
1$ and $m\times 1$ time series, $\bx_t$ is an $r \times 1$ latent
factor process, $\bvarepsilon_t\sim \textrm{WN}(0,{\bf
\Sigma}_\varepsilon)$ is a white noise with zero mean and covariance
matrix ${\bf \Sigma}_\varepsilon$ and $\bvarepsilon_t$ is
uncorrelated with $(\bz_t, \bx_t)$, ${\bf D}$ is an unknown
regression coefficient matrix, and ${\bf A}$ is an unknown factor
loading matrix. The number of the latent factors $r$ is an unknown
(fixed) constant. With the observations $\{ (\by_t, \bz_t): \; t=1,
\ldots, T \}$, the goal is to estimate ${\bf D}, \; {\bf A}$ and
$r$, and to recover the factor process $\bx_t$, when $p$ is large in
relation to the sample size $T$. As our inference will be based on
the serial dependence of each and across $\by_t, \bz_t$ and $\bx_t$,
we assume $E({\bf z}_t)={\bf 0}$ and $E({\bf x}_t)={\bf 0}$ for
simplicity.

In this section, we consider the simple case when $\bz_t$ and
$\bx_t$ are uncorrelated. This condition ensures that the
coefficient matrix $\bD$ in (\ref{eq:nonlinear1}) is identifiable.
However the factor loading matrix $\bA$ and the factor $\bx_t$ are
not uniquely determined by (\ref{eq:nonlinear1}), as we may replace
$({\bf A},{\bf x}_t)$ by $({\bf A}{\bf H},{\bf H}^{-1}{\bf x}_t)$
for any invertible matrix ${\bf H}$. Nevertheless the linear space
spanned by the columns of ${\bf A}$, denoted by $\mathcal{M}({\bf
A})$, is uniquely defined. $\mathcal{M}({\bf A})$ is called the
factor loading space. Hence there is no loss of the generality in
assuming that $\bA$ is a half orthogonal matrix in the sense that
$\bA^\T\bA=\bI_{r}$. In this paper, we always adhere with this
assumption. Once we have specified a particular $\bA$, $\bx_t$ is
uniquely defined accordingly. On the other hand, when
$\textrm{cov}(\bz_t, \bx_t)\ne {\bf 0}$, the endogeneity makes $\bD$
unidentifiable,
% but $\mathcal{M}({\bf A})$ can be still uniquely defined via
% (\ref{eq:nonlinear1})},
which will be dealt with in Section 3 below.

\subsection{Estimation}

Formally the estimation for $\bD$ may be treated as a standard least
squares problem, since
\begin{equation} \label{b2}
{\bf y}_t={\bf D}{\bf z}_t+\beeta_t, \qquad \beeta_t= \bA
\bx_t+\bvarepsilon_t,
\end{equation}
 and cov$(\bz_t, \beeta_t) ={\bf 0}$; see (\ref{eq:nonlinear1}).
Write ${\bf D}=({\bf d}_1,\ldots,{\bf d}_p)^\T$. The least squares
estimator for ${\bf D}$ can be expressed as
\begin{equation} \label{b3}
\widehat{{\bf D}}=(\widehat{{\bf d}}_1,\ldots,\widehat{{\bf
d}}_p)^\T, \qquad \widehat{{\bf d}}_i =
\bigg(\frac{1}{T}\sum_{t=1}^T{\bf z}_t{\bf
z}_t^\T\bigg)^{-1}\bigg(\frac{1}{T}\sum_{t=1}^T y_{i,t}\, {\bf
z}_t\bigg),
\end{equation}
where $y_{i,t}$ is the $i$th component of $\by_t$.

The estimation for $\calM(\bA)$ is based on the residuals $\wh
\beeta_t = \by_t - \wh \bD \bz_t$, using the same idea as
\cite{LamYaoBathia_Biometrika_2011} and \cite{LamYao_AOS_2012},
though we do not assume that the processes concerned are stationary.
To this end, we introduce some notation first. Let
\[
\bSigma_x(k)=\frac{1}{T-k}\sum_{t=1}^{T-k}\textrm{cov}({\bf
x}_{t+k},{\bf x}_t), \qquad
\bSigma_{x\varepsilon}(k)=\frac{1}{T-k}\sum_{t=1}^{T-k}\textrm{cov}({\bf
x}_{t+k},\bvarepsilon_t),
\]
\[
\bSigma_{\eta}(k)=\frac{1}{T-k}\sum_{t=1}^{T-k}\textrm{cov}(\beeta_{t+k},\beeta_t).
\]
When, for example, $\bx_t$ is stationary, $\bSigma_x(k)$ is the
autocovariance matrix of $\bx_t$ at lag $k$. It follows from the
second equation in (\ref{b2}) that for any $k\ne 0$,
\begin{equation}
\bSigma_{\eta}(k)={\bf A}\bSigma_{x}(k){\bf A}^\T+{\bf
A}\bSigma_{x\varepsilon}(k).\label{eq:eta}
\end{equation}
For a prescribed fixed positive integer $\bar{k}$, define
\begin{equation}
{\bf
M}=\sum_{k=1}^{\bar{k}}\bSigma_{\eta}(k)\bSigma_{\eta}(k)^\T.\label{eq:m}
\end{equation}
We assume rank$(\bM)=r$. This is reasonable as it effectively
assumes that the latent factor process $\bx_t$ is genuinely
$r$-dimensional. Since $\bM$ is implicitly sandwiched by $\bA$ and
$\bA^\T$, $\bM \bb={\bf 0}$ for any $\bb \perp \calM(\bA)$. Thus we
may take the eigenvectors of $\bM$ corresponding to non-zero
eigenvalues as the columns of $\bA$, as the choice of $\bA$ is
almost arbitrary as long as $\calM(\bA)$ does not change. Let $\bA
=( \ba_1, \ldots, \ba_{r})$, where $\ba_1, \ldots, \ba_{r}$ be the
$r$ orthonormal eigenvectors of $\bM$ corresponding to the $r$
largest eigenvalues $\lambda_1 \ge \cdots \ge \lambda_{r}>0$. Then
$\bA$ is a half orthogonal matrix in the sense that
$\bA^\T\bA=\bI_{r}$. In the sequel, we always use  $\bA$ defined
this way. When the $r$ non-zero eigenvalues of $\bM$ are distinct,
$\bA$ is unique if we ignore the trivial replacements of $\ba_j$ by
$-\ba_j$.

Let $\wh \bfeta_t =\by_t - \wh \bD \bz_t$ and
\[
\wh \bSigma_\eta (k) =  {1 \over T-k} \sum_{t=1}^{T-k} (\wh
\bfeta_{t+k}-\bar{\bfeta}) (\wh\bfeta_t-\bar{\bfeta})^\T, \qquad
\bar{\bfeta}=\frac{1}{T}\sum_{t=1}^T\wh \bfeta_t.
\]
The above discussion leads to a natural estimator of $\bA$
denoted by  $\wh \bA \equiv
(\wh \ba_1, \ldots, \wh \ba_{r})$. Here $\wh \ba_1, \ldots, \wh
\ba_{r}$ are the orthonormal eigenvectors of $\wh \bM$ corresponding
to the $r$ largest eigenvalues $\wh \lambda_1 \ge \cdots \ge \wh
\lambda_{r}$, where
\begin{equation} \label{eq:hatm}
\wh \bM = \sum_{k=1}^{\bar k} \wh \bSigma_\eta (k) \wh \bSigma_\eta
(k)^\T.
% \; {\rm for} \; \wh \bSigma_\eta (k) =  {1 \over T-k}
% \sum_{t=1}^{T-k} (\wh \bfeta_{t+k}-\bar{\bfeta})
% (\wh\bfeta_t-\bar{\bfeta})^\T\;\textrm{and}\;\bar{\bfeta}=\frac{1}{T}\sum_{t=1}^T\wh
% \bfeta_t.
\end{equation}
 Since $\wh \bA$ is a half orthogonal matrix, we may extract the
factor process by $\wh \bx_t = % \wh {\bf A}^\T \wh \bfeta_t =
\wh {\bf A}^\T(\by_t - \wh \bD \bz_t)$; see (\ref{b2}).

All the arguments above are based on
a known $r$ which is actually unknown in practice. The determination of $r$ is
a key step in our inference.
In practice we may
estimate it by the ratio estimator
\begin{equation}
\widehat{r}=\arg\min
\bigg\{\frac{\widehat{\lambda}_{j+1}}{\widehat{\lambda}_{j}} :1\leq
j\leq R \bigg\},\label{eq:number}
\end{equation}
where  $\wh \lambda_1\geq \cdots\geq \wh \lambda_p$ are the
eigenvalues of $\wh \bM$, and $R$ is a constant which may be taken
as $R=p/2$; see \cite{LamYao_AOS_2012} for further discussion on
this estimation method.

\subsection{Asymptotic properties}

We present the asymptotic theory for the estimation methods
described in Section 2.2 above when $T, \; p \to \infty$ while $r$
is fixed. We also assume $m$ fixed now; see Section 4 below for the
results when $m\to \infty$ as well. We do not impose stationarity
conditions on $\by_t, \bz_t $ and $\bx_t$. Instead we assume that
they are mixing processes; see Condition~\ref{as1} below.
 Hence  our results in the special case when $\bz_t \equiv \bf0$
% for the estimation of $\bA$ and $r$
extend those in \cite{LamYaoBathia_Biometrika_2011} and
\cite{LamYao_AOS_2012} to nonstationary cases.
\cite{PanYao_Bioka_2008} dealt with a different method for
nonstationary factor models.

We introduce some notation first. For any matrix $\bH$, we denote by
$\|\bH\|_F = \{{\rm tr}(\bH^\T\bH)\}^{1/2}$ the Frobenius norm of
$\bH$, and by $\|\bH\|_2 = \{\lambda_{\max}(\bH^\T\bH)\}^{1/2}$ the
$L_2$-norm,  where tr$(\cdot)$ and $\lambda_{\max}(\cdot)$ denote,
respectively, the trace and the maximum eigenvalue of a square
matrix. We also %write $\lambda_{\max} (\bH) = \|\bH\|_2$ and
denote by $\|\bH\|_{\min}$ the square-root of the minimum  nonzero
eigenvalue of $\bH^\T\bH$. Note that when $\bH= \bh$ is a vector,
$\|\bh\|_F= \|\bh\|_2 = \|\bh\|_{\min}= (\bh^\T\bh)^{1/2}$, i.e. the
conventional Euclidean norm for vector $\bh$.

\begin{as}
\label{as1} The process $\{ (\by_t, \bz_t, \bx_t) \}$ is
$\alpha$-mixing with the mixing coefficients satisfying the
condition $ \sum_{k=1}^\infty\alpha(k)^{1-2/\gamma}<\infty$ for
some $\gamma>2$, where
$$
\alpha(k) = \sup_i  \sup_{ A \in \calF_{-\infty}^i, \; B \in
\calF_{i+k}^\infty} \big| P(A \cap B) - P(A) P(B) \big|,
$$
and $\calF_i^j$ is the $\sigma$-field generated by $\{ (\by_t,
\bz_t, \bx_t): \; i \le t \le j \}$.
\end{as}

\begin{as}
\label{as2} For any $i=1,\ldots,m$, $j=1,\ldots,p$ and $t$,
 $
E(|z_{i,t}|^{2\gamma})\leq C_1$, $E(|\zeta_{j,t}|^{2\gamma})\leq
C_1$ and $E(|\varepsilon_{j,t}|^{2\gamma})\leq C_1$, where $C_1>0$
is a constant, $\gamma$ is given in Condition \ref{as1}, and
$z_{i,t}$ is the $i$th element of ${\bf z}_t$, $\zeta_{j,t}$ and
$\varepsilon_{j,t}$ are the $j$th element of, respectively,
$\bA\bx_t$ and $\bve_t$.
\end{as}

\begin{as}
\label{as3}  There exists a constant $C_2 >0$ such that
 $\lambda_{\textrm{min}}\{E({\bf z}_t{\bf z}_t^\T)\}>C_2$ for all $t$.
\end{as}

Condition $E(|\zeta_{j,t}|^{2\gamma})\leq
C_1$ in Condition
 \ref{as2} can be guaranteed by some suitable conditions on each
$x_{i,t}$, as  ${\bf A}$ is a half orthogonal matrix. For example, it holds if
$\max_{i,t} E(|x_{i,t}|^{2\gamma})<\infty$.
Proposition~\ref{pn1} below
establishes the convergence rate of the estimator for the $p\times m$
coefficient matrix $\bD$.  Since $p \to \infty$ together with the sample
size $T$, the convergence rate depends on~$p$. Especially when $p/T \to
0$, the least squares estimator $\wh \bD$ is a consistent estimator for
$\bD$. This condition can be relaxed if we impose some sparse condition
on $\bD$, and then apply appropriate thresholding on $\wh \bD$. We do not
pursue this further here. When $p$ is fixed, the convergence rate is
$T^{1/2}$ which is the optimal rate for the regression with the dimension fixed.

\begin{pn}
\label{pn1} Let Conditions \ref{as1}-\ref{as3} hold. As $T \to
\infty$ and $ p \to \infty$, it holds that
\[
\|\widehat{{\bf D}}-{\bf D}\|_F=O_p\big(p^{1/2}T^{-1/2}\big).
\]
\end{pn}

To state the results for estimating factor loadings, we introduce
more conditions.

%\begin{as}
%\label{as4} There exists $0\leq \delta\leq1$ such that $\|{\bf
%a}_i\|_2^2\asymp p^{1-\delta}$ and $\min_{\theta_j,j\neq i}\|{\bf
%a}_i-\sum_{j\neq i}\theta_j{\bf a}_j\|_2^2\asymp p^{1-\delta}$ for
%all $i=1,\ldots,r$.
%\end{as}

\begin{as}
\label{as5} There exist positive constants $C_i~(i=3,4)$ and
$\delta\in[0,1]$ such that $
C_3p^{1-\delta}\leq\|\bSigma_x(k)\|_{\min}\leq\|\bSigma_x(k)\|_2\leq
C_4p^{1-\delta}$ for all
% C_3\leq\|\bSigma_x(k)\|_{\min}\leq\|\bSigma_x(k)\|_2\leq C_4$ for all
$k=1,\ldots,\bar{k}$.
\end{as}

\begin{as}
\label{as6} Matrix $\bM$ admits $r$ distinct positive eigenvalues
$\lambda_1 > \cdots > \lambda_{r}>0$.
\end{as}

The constant $\delta$ in Condition \ref{as5} controls the strength of the
factors. When $\delta=0$, the factors are strong. When $\delta >0$, the
factors are weak. In fact the value of $\delta$ reflects the sparse level
of the factor loading matrix $\bA$, and a certain degree of sparsity is
present when $\delta >0$. Therefore not all components of $\by_t -
\bD{\bf z}_t$ carry the information for all factor components. This
causes difficulties in recovering the factor process. This argument will
be verified in Theorem 2.2. See also Remark 1 in \cite{LamYao_AOS_2012}.
Condition \ref{as6} implies that $\bA$ defined as in Section 2.2 above is unique.
This simplifies the presentation significantly,
as Theorem~\ref{tm1} below can present the convergence rates of
the estimator for $\bA$ directly. Without condition \ref{as6}, the same convergence
rates can be obtained for the estimation of the linear space
$\mathcal{M}(\widehat{\bA})$; see (\ref{ratesSP}) below.
Let
\[
\kappa_{1}=\min_{1\leq
k\leq\bar{k}}\|\bSigma_{x\varepsilon}(k)\|_{\min}~~\textrm{and}~~\kappa_2=\max_{1\leq
k\leq\bar{k}}\|\bSigma_{x\varepsilon}(k)\|_2.
\]
Note that both $\kappa_{1}$ and $\kappa_2$ may diverge as $p\to
\infty$.

\begin{tm}
\label{tm1} Let Conditions \ref{as1}-\ref{as6} hold. Suppose that
$r$ is known and fixed, then
\[
\|\widehat{{\bf A}}-{\bf A}\|_2 =~\left\{
     \begin{array}{ll}
       O_p(p^\delta
T^{-1/2}), & {\rm if\; }  \kappa_2=o(p^{1-\delta}) {\rm \;\; and
\;\;}
p^{2\delta} T^{-1}=o(1); \\[1ex]
       O_p(\kappa_1^{-2}\kappa_2pT^{-1/2}), \quad &{\rm if\;} p^{1-\delta}=o(\kappa_1)
{\rm \;\; and \;\;} \kappa_1^{-2}\kappa_2pT^{-1/2}=o(1).
     \end{array}
   \right.
\]
\end{tm}

The convergence rates in Theorem~\ref{tm1} above are exactly the
same as Theorem~1 of \cite{LamYaoBathia_Biometrika_2011} which deals
with a pure factor model, i.e. model (\ref{b2}) with $\bz_t\equiv
\bf0$. In this sense, the estimator $\wh \bA$ is asymptotically
adaptive to unknown $\bD$.

\begin{tm}
\label{tm2} Let Conditions \ref{as1}-\ref{as6} hold, and $r$ be
known and fixed. If $\|\bSigma_{\varepsilon}\|_2$ is bounded as
$p\to \infty$, then
\[
p^{-1/2}\|\widehat{{\bf A}}\widehat{{\bf x}}_t-{\bf A}{\bf
x}_t\|_2=O_p(\|\widehat{{\bf A}}-{\bf A}\|_2+p^{-1/2}+T^{-1/2}).
\]
\end{tm}

Theorem~\ref{tm2} deals with the convergence of the extracted factor
term. Combining it with Theorem \ref{tm1}, we obtain
\begin{align*}
& p^{-1/2}\|\widehat{{\bf A}}\widehat{{\bf x}}_t-{\bf A}{\bf
x}_t\|_2\\[1ex]
=&\left\{
     \begin{array}{ll}
       O_p(p^\delta T^{-1/2}+p^{-1/2}),&{\rm if}\;
\kappa_2=o(p^{1-\delta}) {\rm \; \; and \;\;}  p^{2\delta}T^{-1}=o(1); \\[1ex]
       O_p(\kappa_1^{-2}\kappa_2pT^{-1/2}+p^{-1/2}+T^{-1/2}), \quad&{ \rm if} \;
 p^{1-\delta}=o(\kappa_1) {\rm \; \; and \;\;} \kappa_1^{-2}\kappa_2pT^{-1/2}=o(1).
     \end{array}
   \right.
\end{align*}
Thus when all the factors are strong (i.e.  $\delta=0$) and
$\kappa_2=o(p)$, it holds that $p^{-1/2}\|\widehat{{\bf
A}}\widehat{{\bf x}}_t-{\bf A}{\bf x}_t\|_2
=O_p(p^{-1/2}+T^{-1/2})$, which is the optimal convergence rate
specified in Theorem 3 of \cite{Bai_Econometrica_2003}.

In general the choice of $\bA$ in model (\ref{eq:nonlinear1}) is not
unique, we consider the error in estimating $\calM(\bA)$ instead of
a particular $\bA$, as $\calM(\bA)$ is uniquely defined by
(\ref{eq:nonlinear1}) and does not vary with different choices of
$\bA$. To this end, we adopt the discrepancy measure used by
\cite{PanYao_Bioka_2008}: for two $p\times r$ half orthogonal
matrices ${\bf H_1}$ and ${\bf H}_2$ satisfying the condition ${\bf
H}_1^\T{\bf H}_1={\bf H}_2^\T{\bf H}_2=\bI_{r}$, the difference
between the two linear spaces $\mathcal{M}({\bf H}_1)$ and
$\mathcal{M}({\bf H}_2)$ is measured by
\begin{equation}
D(\mathcal{M}({\bf H}_1),\mathcal{M}({\bf
H}_2))=\sqrt{1-\frac{1}{r}\textrm{tr}({\bf H}_1{\bf H}_1^\T{\bf
H}_2{\bf H}_2^\T)}.\label{eq:D}
\end{equation}
In fact $D(\mathcal{M}({\bf H}_1),\mathcal{M}({\bf H}_2))$ always
takes values between 0 and 1. It is equal to $0$ if and only if
$\mathcal{M}({\bf H}_1)=\mathcal{M}({\bf H}_2)$, and to $1$ if and
only if $\mathcal{M}({\bf H}_1)\perp \mathcal{M}({\bf H}_2)$.
% \cite{PanYao_Bioka_2008} shows that this discrepancy is a metric
% equipped on a quotient space where treat ${\bf H}_1$ and ${\bf H}_2$
% are the same if ${\bf H}_1={\bf H}_2{\bf V}$ for some orthogonal
% matrix ${\bf V}$. For the general linear space spanned by some
% orthogonal functions, the similar definition of discrepancy between
% two linear spaces can be found in
% \cite{BathiaYaoZiegelmann_AOS_2010}. The following theorem gives the
% bounds for the discrepancy between $\mathcal{M}(\widehat{{\bf Q}})$
% and $\mathcal{M}({\bf Q})$.

\begin{tm}
\label{tm3} Let Conditions \ref{as5}-\ref{as6} hold. Suppose that
$r$ is known and fixed, then
\[
\{ D(\mathcal{M}(\widehat{\bf A}),\mathcal{M}({\bf A}))\}^2
\asymp\|(\widehat{{\bf A}}-{\bf A})^\T(\widehat{{\bf A}}-{\bf
A})-{\bf A}^\T(\widehat{{\bf A}}-{\bf A})(\widehat{{\bf A}}-{\bf
A})^\T{\bf A}\|_2.
\]
\end{tm}

This theorem establishes the link between
$D(\mathcal{M}(\widehat{\bf A}),\mathcal{M}({\bf A}))$ and
$\widehat{{\bf A}}-{\bf A}$ when $r$ is known.  Obviously, the
%quantity on the right hand side
RHS of the above expression can be bounded by $2\|\widehat{{\bf A}}-{\bf
A}\|_2^2$. This implies that $D(\mathcal{M}(\widehat{\bf
A}),\mathcal{M}({\bf A}))=O_p(\|\widehat{{\bf A}}-{\bf A}\|_2)$.
In fact, the convergence of $D(\mathcal{M}(\widehat{\bf
A}),\mathcal{M}({\bf A}))$ does not depend on Condition 2.5.
Even when ${\bf M}$ admits multiple non-zero eigenvalues, and,
therefore, ${\bf A}$ is not uniquely defined,
it can be shown based on the similar arguments as for Theorem 1 in
\cite{ChangGuoYao_2014}
that
\begin{equation} \label{ratesSP}
D(\mathcal{M}(\widehat{\bf A}),\mathcal{M}({\bf A})) =~\left\{
     \begin{array}{ll}
       O_p(p^\delta
T^{-1/2}), & {\rm if\; }  \kappa_2=o(p^{1-\delta}) {\rm \;\; and \;\;}
p^{2\delta} T^{-1}=o(1); \\[1ex]
       O_p(\kappa_1^{-2}\kappa_2pT^{-1/2}), \quad &{\rm if\;} p^{1-\delta}=o(\kappa_1)
{\rm \;\; and \;\;} \kappa_1^{-2}\kappa_2pT^{-1/2}=o(1),
     \end{array}
   \right.
\end{equation}
which is the same as that followed by Theorem \ref{tm3} when Condition 2.5 holds.

 Theorems \ref{tm1}-\ref{tm3} above present the asymptotic properties
when the number of factors $r$ is assumed to be known. However, in
practice we need to estimate $r$ as well. \cite{LamYao_AOS_2012}
showed that for the ratio estimator  $\wh r$ defined in
(\ref{eq:number}), $P(\wh r \ge r) \to 1.$ In spite of favorable
finite sample evidences reported in \cite{LamYao_AOS_2012}, it
remains as a unsolved challenge to establish the consistency $\wh
r$. Following the idea of \cite{Xiaetal_2013}, we adjust the ratio
estimator as follows
\begin{equation} \label{b9}
\wt r = \arg\min \bigg\{\frac{\widehat{\lambda}_{j+1} +
C_T}{\widehat{\lambda}_{j} + C_T} :1\leq j\leq R \bigg\},
\end{equation}
where $C_T = (p^{1-\delta}+\kappa_2)pT^{-1/2}\log T$.
Theorem~\ref{tm4} shows that $\wt r$ is a consistent estimator for
$r$.

\begin{tm}
\label{tm4} Let Conditions \ref{as1}-\ref{as6} hold, and $
(p^{1-\delta}+\kappa_2)pT^{-1/2}\log T=o(1)$. Then $P(\wt r \ne r)
\to 0$.
\end{tm}

With the estimator $\wt r$, we may define an estimator for $\bA$ as
$\wt \bA = (\wh \ba_1, \ldots, \wh \ba_{\wt r} )$, where $\wh \ba_1,
\ldots, \wh \ba_{\wt r}$ are the orthonormal eigenvectors of $\wh
\bM$, defined in (\ref{eq:hatm}), corresponding to the $\wt r$
largest eigenvalues. Then $\wt \bA = \wh \bA$ when $\wt r = r$. To
measure the error in estimating the factor loading space, we use
\[
\widetilde{D}(\mathcal{M}(\widetilde{{\bf A}}),\mathcal{M}({\bf
A}))=\sqrt{1-\frac{1}{\max(\widetilde{r},r)}\textrm{tr}(\widetilde{{\bf
A}}\widetilde{{\bf A}}^\T{\bf A}{\bf A}^\T)}.
\]
This is a modified version of (\ref{eq:D}). It takes into account
the fact that the dimensions of $\calM(\wt \bA)$ and $\calM(\bA)$ may
be different. Obviously
 $\widetilde{D}(\mathcal{M}(\widetilde{{\bf
A}}),\mathcal{M}({\bf A}))=D(\mathcal{M}(\widehat{{\bf
A}}),\mathcal{M}({\bf A}))$ if $\widetilde{r}=r$. We show below that
$\widetilde{D}(\mathcal{M}(\widetilde{{\bf A}}),\mathcal{M}({\bf
A})) \to 0 $ in probability at the same rate as
$D(\mathcal{M}(\widehat{{\bf A}}),\mathcal{M}({\bf A}))$. Hence even
without knowing $r$, $\calM(\wt \bA)$ is a consistent estimator for
$\calM(\bA)$. Let $\rho = \rho(T,p)$ denote the convergence rate of
$D(\mathcal{M}(\widehat{{\bf A}}),\mathcal{M}({\bf A}))$, i.e. $\rho
D(\mathcal{M}(\widehat{{\bf A}}),\mathcal{M}({\bf A})) =O_p(1),$ see
Theorems \ref{tm1} and \ref{tm3}. For any $\epsilon>0$, there exists
a positive constant $M_{\epsilon}$ such that
$P\{\rho{D}(\mathcal{M}(\widehat{{\bf A}}),\mathcal{M}({\bf
A}))>M_{\epsilon}\}<\epsilon$. Then,
\[
\begin{split}
&~P\{\rho\widetilde{D}(\mathcal{M}(\widetilde{{\bf
A}}),\mathcal{M}({\bf A}))>M_{\epsilon}\}\\
\leq&~P\{\rho{D}(\mathcal{M}(\widehat{{\bf A}}),\mathcal{M}({\bf
A}))>M_{\epsilon},\widetilde{r}=r\}+P\{\rho
\widetilde{D}(\mathcal{M}(\widetilde{{\bf
A}}),\mathcal{M}({\bf A}))>M_{\epsilon},\widetilde{r}\neq r\}\\
\le &~P\{\rho D(\mathcal{M}(\widehat{{\bf A}}),\mathcal{M}({\bf
A}))>M_{\epsilon}\}+o(1)\\
\le&~ \epsilon+o(1) \to \epsilon\\
\end{split}
\]
which implies $\rho\widetilde{D}(\mathcal{M}(\widetilde{{\bf
A}}),\mathcal{M}({\bf A}))=O_p(1)$. Hence,
$\widetilde{D}(\mathcal{M}(\widetilde{{\bf A}}),\mathcal{M}({\bf
A})) \to 0 $ shares the same convergence rate of
$D(\mathcal{M}(\widehat{{\bf A}}),\mathcal{M}({\bf A}))$ which means
that $\mathcal{M}(\widetilde{{\bf A}})$ has the oracle property in
estimating the factor loading space $\mathcal{M}({\bf A})$.

\section{Models with endogeneity}

In last section, the consistent estimation for the coefficient
matrix $\bD$ is used in identifying the latent factor process. The
consistency is guaranteed by the assumption that cov$(\bz_t, \bx_t)
= E(\bz_t \bx_t^\T)={\bf 0}$. However when the endogeneity exists in
model (\ref{eq:nonlinear1}) in the sense that the regressor $\bz_t$
and the latent factor $\bx_t$ are contemporaneously correlated with
each other, $\bD$ is no longer identifiable. Nevertheless
(\ref{eq:nonlinear1}) can be written as
\begin{align} \label{d1}
\by_t &\; = \; [\bD+ \bA E(\bx_t \bz_t^\T) \{E(\bz_t
\bz_t^\T)\}^{-1}]\bz_t +
\bA[\bx_t - E(\bx_t \bz_t^\T) \{E(\bz_t \bz_t^\T)\}^{-1}\bz_t] + \bve_t \\
&\; \equiv \; \bD^\star \bz_t + \bA \bx^\star_t + \bve_t, \nonumber
\end{align}
where the latent factor $\bx^\star_t = \bx_t - E(\bx_t \bz_t^\T)
\{E(\bz_t \bz_t^\T)\}^{-1}\bz_t$ is uncorrelated with the regressor
$\bz_t$. Hence if we apply the methods presented in Section 2 to
model (\ref{eq:nonlinear1}) in the presence of the endogeneity, $\wh
\bD$ defined in (\ref{b3}) is a consistent estimator for $\bD^\star
= \bD+ \bA E(\bx_t \bz_t^\T) \{E(\bz_t \bz_t^\T)\}^{-1}$ instead of
the original regression coefficient $\bD$, provided that $\bD^\star$
so defined is a constant matrix independent of $t$. The latter is
guaranteed when both $\bx_t$ and $\bz_t$ are stationary.
Furthermore, the recovered factor process $\wh \bx_t$ is an
estimator for $\bx^\star_t$.
% the asymptotic properties presented in Section~2.3 still hold (with
% some obvious alternations on the corresponding regularity
% conditions),
Hence in the presence of the endogeneity and if $\bD^\star$ defined in (\ref{d1})
is a constant matrix,  the factor
loading space $\calM(\bA)$ can still be estimated consistently
although the ordinary least squares estimator for the regression coefficient
matrix $\bD$ is no longer consistent.
% {\bf From (\ref{d1}), we can find that although the
%consistent estimation for ${\bf D}$ does not hold if we ignore the
%endogeneity in the estimation procedure, we still can obtain the
%consistent estimation for $\bA$ and $\mathcal{M}(\bA)$.}

For some applications, the interest lies in estimating the
`original' $\bD$ and $\bx_t$; see, e.g., \cite{Angrist_1991}. Then
we may employ a set of instrument variables ${\bf w}_t$ in the sense
that ${\bf w}_t$ is correlated with $\bz_t$ but uncorrelated with
both $\bx_t$ and $\bve_t$. Usually, we require that $\bw_t$ is
$q\times 1 $ with $q\ge m$. It follows from (\ref{eq:nonlinear1})
that
\begin{equation} \label{d2}
\by_t \bw_t^\T = \bD \bz_t \bw_t^\T + \bve^\star_t, \quad
\bve^\star_t = \bA \bx_t \bw_t^\T + \bve_t\bw_t^\T.
\end{equation}
Since $E(\bx_t \bw_t^\T)={\bf 0}$ and $E(\bve_t\bw_t^\T)={\bf 0}$,
we may view the first equation in the above expression as similar to
a `normal equation' in a least squares problem by ignoring
$\bve_t^\star$. This leads to the following estimator for $\bD$:
\begin{equation}
\widehat{{\bf D}}=\bigg(\frac{1}{T}\sum_{t=1}^\T \by_t \bw_t^\T {\bf
R}^\T\bigg)\bigg(\frac{1}{T}\sum_{t=1}^T{\bf z}_t{\bf w}_t^\T{\bf
R}^\T\bigg)^{-1}.\label{eq:ive}
\end{equation}
where ${\bf R}$ is any $m\times q$ constant matrix with rank$({\bf
R}) = m$,
 to match the lengths of
$\bw_t$ and $\bz_t$. When $q=m$, we can choose ${\bf R}= \bI_m$.
This is the `instrument variables method' widely used in
econometrics. We refer to \cite{Morimune_1983},
\cite{BoundJaegerBaker_1996}, \cite{DonaldNewey_2001},
\cite{HahnHausman_2002} and \cite{CanerFan_2012} for further
discussion on the choice of instrument variables and the related
issues. It follows from (\ref{d2}) and (\ref{eq:ive}) that
\[
\widehat{{\bf D}} - \bD = \bigg(\frac{1}{T}\sum_{t=1}^T
\bve_t^\star{\bf
 R}^\T\bigg)\bigg(\frac{1}{T}\sum_{t=1}^T{\bf z}_t{\bf
w}_t^\T{\bf R}^\T\bigg)^{-1}.
\]
The proposition below shows that $\wh \bD$ is a consistent estimator
with the optimal convergence rate. See also Proposition \ref{pn1}.

\begin{as}
\label{as7} For any $i=1,\ldots,q$ and $t$,
$E(|w_{i,t}|^{2\gamma})\leq C_1$ for  $\gamma>2$ and $C_1>0$
specified in, respectively, Conditions \ref{as1} and \ref{as2}.
\end{as}

\begin{as}
\label{as8} The smallest eigenvalue of $\{E({\bf w}_t{\bf
z}_t^\T)\}^\T{\bf R}^\T{\bf R}\{E({\bf w}_t{\bf z}_t^\T)\}$ is
uniformly bounded away from zero for all $t$.
\end{as}

% Condition \ref{as7} requires each component of ${\bf w}_t$ has
% $2\gamma$ moments.
Condition \ref{as8} implies that all the components of the
instrument variables $\bw_t$ are correlated with the regressor
$\bz_t$. When $q=m$ and $\bR=\bI_m$, it reduces to the condition
that all the singular values of $E({\bf w}_t{\bf z}_t^\T)$ are
uniformly bounded away from zero for all $t$.
%  which implies ${\bf w}_t$ should be correlated to ${\bf z}_t$. When
% $j>m$, Condition \ref{as8} means the new structured instrument
% variables ${\bf R}{\bf w}_t$ are correlated to ${\bf z}_t$.

\begin{pn}
\label{pn2} Let Conditions \ref{as1}-\ref{as2} and
\ref{as7}-\ref{as8} hold. As $T\rightarrow\infty$ and
$p\rightarrow\infty$, it holds that
\[
\|\widehat{{\bf D}}-{\bf D}\|_F=O_p(p^{1/2}T^{-1/2}).
\]
\end{pn}

With the consistent estimator $\wh \bD$ in (\ref{eq:ive}), the
factor loading space and the latent factor process may be estimated
in the same manner as in Section 2.2. The asymptotic properties
presented in Theorems \ref{tm1}-\ref{tm3} can be reproduced in the
similar manner.

\section{Models with nonlinear regression functions}\label{sn}

Now we consider the model with nonlinear regression term:
\begin{equation}
\by_t = \bg(\bu_t) + \bA \bx_t + \bve_t,
        \label{eq:nonlinear}
\end{equation}
where ${\bf g}(\cdot)$ is an unknown nonlinear function, ${\bf u}_t$
is an observed process with fixed dimension, and other terms are the
same as in model (\ref{eq:nonlinear1}). One way to handle a
nonlinear regression is to transform it into a high-dimensional
linear regression problem. To this end, let $\bg = (g_1, \ldots,
g_p)^\T$, and
\[
g_i(\bu) = \sum_{j=1}^\infty d_{i,j} l_j(\bu), \qquad i =1, 2,
\ldots,
\]
where $\{l_j(\cdot) \}$ is a set of base functions. Suppose we use
the approximation with the first $m$ terms only. Let ${\bf
z}_t=(l_1({\bf u}_t),\ldots,l_m({\bf u}_t))^\T$, and $\bD$ be the
$p\times m$ matrix with $d_{i,j}$ as its $(i,j)$-th element, then
(\ref{eq:nonlinear}) can be expressed as
\begin{equation} \label{eq:newnonlinear}
\by_t = \bD \bz_t + \bA \bx_t + \bve_t + \bfe_t,
\end{equation}
where the additional error term $\bfe_t$ collects the residuals in
approximating $\bg(\cdot)$ by the first $m$ terms only, i.e. the
$i$th component of $\bfe_t$ is $\sum_{j>m} d_{i,j} l_j(\bu_t)$. This
makes (\ref{eq:newnonlinear}) formally different from model
 (\ref{eq:nonlinear1}).
Furthermore a fundamentally new feature in (\ref{eq:newnonlinear})
is that $m$ may be large in relation to $p$ or/and $T$. Hence the
new asymptotic theory with all $T, p, m \to \infty$ together will be
established  in order to take into account those non-trivial
changes. Due to (\ref{d1}), we may always assume that
% {\bf the identification
% condition $E(\bA {\bf x}_t|{\bf u}_t)={\bf 0}$ holds for each $t$,
% which implies}
cov$(\bz_t, \bx_t)={\bf 0}$. Condition~\ref{as13} below ensures that
$\bfe_t$ in (\ref{eq:newnonlinear}) is asymptotically negligible.
Hence model (\ref{eq:newnonlinear}) is as identifiable as
(\ref{eq:nonlinear1}) at least asymptotically when $m \to \infty$.
Consequently we may estimate $\bD$ using  the ordinary least squares
estimator:
\[
\widehat{{\bf D}}=\bigg(\frac{1}{T}\sum_{t=1}^T{\bf y}_t{\bf
z}_t^\T\bigg)\bigg(\frac{1}{T}\sum_{t=1}^T{\bf z}_t{\bf
z}_t^\T\bigg)^{-1}.
\]

We introduce some regularity conditions first.
\begin{as}
\label{as12} Supports of the process ${\bf u}_t$ are subsets of
$\mathcal{U}$, where $\mathcal{U}$ is compact with nonempty
interior. Furthermore the density function of ${\bf u}_t$ is
uniformly bounded and bounded away from zero for all $t$.
\end{as}

\begin{as}
\label{as13}  It holds for all large $m$ that
% For given basis functions $\{l_j({\bf
% u})\}_{j=1}^\infty$ and $i=1,\ldots,p$, there exists
% $(d_{i,1},\ldots,d_{i,m})^\T\in\mathbb{R}^m$ such that
\[
\sup_i\sup_{{\bf u}\in\mathcal {U}}\bigg|g_i({\bf
u})-\sum_{j=1}^md_{i,j}l_j({\bf u})\bigg|=O(m^{-\lambda})
\]
where $\lambda>1/2$ is a constant.
\end{as}

\begin{as}
\label{as9} The eigenvalues of $E({\bf z}_t{\bf z}_t^\T)$, are
uniformly bounded away from zero and infinity for all $t$, where
${\bf z}_t=(l_1({\bf u}_t),\ldots,l_m({\bf u}_t))^\T$.
\end{as}

\begin{as}
\label{as10} $E(\bA {\bf x}_t|{\bf u}_t)={\bf 0}$ and
$E(\bvarepsilon_t|{\bf u}_t)={\bf 0}$ for all $t$.
\end{as}

\begin{as}
\label{as11} For each $j=1,\ldots,m$, % the basis function $l_j({\bf u})$ satisfies
$E(|l_j({\bf u}_t)|^{2\gamma})\leq C_1$, where $\gamma>2$ and
$C_1>0$ are specified in, respectively, Conditions \ref{as1} and
\ref{as2}.
\end{as}

Condition \ref{as12} is often assumed in nonparametric estimation,
it can be weakened at the cost of lengthier proofs. Condition
\ref{as13} quantifies the approximation error for regression
function $\bg(\cdot)$. It is fulfilled by
commonly used sieve basis functions such as spline, % power series,
wavelets, or the Fourier series, provided that all components of
$\bg(\cdot)$ are in the H\"{o}lder space. See
\cite{AiChen_Econometrica_2003} for further detail on the sieve method.

\begin{pn}
\label{pnn5} Let Conditions \ref{as1}-\ref{as2} and
\ref{as13}-\ref{as11} hold, and $mT^{-1/2}=o(1)$. Then
\[
\|\widehat{{\bf D}}-{\bf
D}\|_F=O_p(p^{1/2}m^{1/2}T^{-1/2}+p^{1/2}m^{1/2-\lambda}).
\]
\end{pn}

Comparing this proposition with Propositions~\ref{pn1} and
\ref{pn2}, $m$ enters the convergence rates, and the term
$O_p(p^{1/2}m^{1/2-\lambda})$ is due to approximating ${\bf g}({\bf
u}_t)$ by ${\bf D}{\bf z}_t$.
% to the corresponding result for the case where $m$ is
%fixed and $\textrm{cov}({\bf z}_t,{\bf x}_t)={\bf O}$, there is an
%additional term $O_p(p^{1/2}m^{-\lambda})$ in above result. This
%term is generated by the approximation of ${\bf g}({\bf u}_t)$ by
%${\bf D}{\bf z}_t$ which is the bias term in nonparametric estimation.
Based on the estimator $\wh \bD$, we can define an estimator for the
nonlinear regression function
\[
\widehat{{\bf g}}({\bf u})=\widehat{{\bf D}}(l_1({\bf
u}),\ldots,l_m({\bf u}))^\T.
\]
The theorem below follows from Proposition~\ref{pnn5}. It gives the
convergence rate for $\widehat{{\bf g}}$.
\begin{tm}
\label{tm5} Let Conditions \ref{as1}-\ref{as2} and
\ref{as12}-\ref{as11} hold, and $mT^{-1/2}=o(1)$. Then
\[
\int_{{\bf u}\in\mathcal{U}}\|\widehat{{\bf g}}({\bf u})-{\bf
g}({\bf u})\|_2^2 ~d{\bf u}=O_p(pmT^{-1}+pm^{-2\lambda}).
\]
\end{tm}

It is easy to see from Theorem~\ref{tm5} that the best rate for $\wh
\bg(\cdot)$ is attained if we choose $m \asymp T^{1/(2\lambda+1)}$,
which fulfills the condition $mT^{-1/2}=o(1)$ as $\lambda>1/2$. When
$\bg(\cdot)$ is twice differentiable, $\lambda=2$ for some basis
functions, the convergence rate is $pT^{-4/5}$. This is the optimal
rate for the nonparametric regression of $p$ functions
\citep{Stone_1985}. Hereafter, we always set $m \asymp
T^{1/(2\lambda+1)}$.
% $ as
% \[
% m\asymp T^{1/(2\lambda+1)}.
% \]

% As the procedure proposed in last section, we let $
% \widehat{\beeta}_t={\bf y}_t-\widehat{{\bf D}}{\bf z}_t $ and
% $\widehat{\Sigma}_\eta(k)=(T-k)^{-1}\sum_{t=1}^{T-k}\widehat{\beeta}_{t+k}\widehat{\beeta}_t^\T$.
% Define $\widehat{{\bf M}}$ like (\ref{eq:m}). There is a remarkable
% thing is that the selection of $\bar{k}$ in this situation is longer
% a fixed constant, which should be diverging at some rate to obtain
% the optimal convergence rate of $\|\widehat{{\bf Q}}-{\bf Q}\|_2$.
% The following Theorem gives the general result for the convergence
% rate of $\|\widehat{{\bf Q}}-{\bf Q}\|_2$ and the selection of the
% parameter $\bar{k}$ is proposed below this theorem.

With the estimator $\wh \bD$, we may proceed as in Section~2.2 to
estimate the factor loading space and to recover the latent factor
process. However there is a distinctive new feature now: the number
of lags $\bar k$ used in defining both $\bM$ in (\ref{eq:m}) and
$\wh \bM$ in (\ref{eq:number}) may tend to infinity together with
$m$ in order to achieve good convergence rates.
% This point will be further explored after Theorem~\ref{tmn5}.

\begin{tm}
\label{tmn5} Let conditions \ref{as1}-\ref{as2}, \ref{as5} and
\ref{as13}-\ref{as11} hold, $\lambda\geq1$, $\bar{k}T^{-1/2}=o(1)$,
and $m \asymp T^{1/(2\lambda+1)}$.
 Suppose that $r$ is
known, and the $r$ positive eigenvalues of ${\bf M}$ are distinct.
Then
\[
\|\widehat{{\bf A}}-{\bf A}\|_2=\left\{
                                  \begin{aligned}
                                    &O_p\{p^\delta[\bar{k}^{1/2}T^{-1/2}+\bar{k}^{-1}T^{(1-\lambda)/(2\lambda+1)}]\}, \\
                                    &~~~~ \hbox{if $\kappa_2=o(p^{1-\delta})$ and $p^{2\delta}[\bar{k}T^{-1}+T^{(2-2\lambda)/(2\lambda+1)}]=o(1)$;} \\
                                    &O_p\{p\kappa_2\kappa_1^{-2}[\bar{k}^{1/2}T^{-1/2}+\bar{k}^{-1}T^{(1-\lambda)/(2\lambda+1)}]\},\\
                                    &~~~~ \hbox{if $p^{1-\delta}=o(\kappa_1)$ and $p^{2}\kappa_2^2\kappa_1^{-4}[\bar{k}T^{-1}+T^{(2-2\lambda)/(2\lambda+1)}]=o(1)$.}
                                  \end{aligned}
                                \right.
\]
\end{tm}

From Theorem \ref{tmn5}, the best convergence rate for $\wh \bA$ is
attained when we choose $\bar{k}\asymp T^{1/(2\lambda+1)}$. The
model with linear regression considered in Section 2.3 corresponds
to the cases with $\lambda=\infty$. Note Theorem \ref{tmn5} implies
that $\bar{k}\asymp 1$ should be used when $\lambda=\infty$ and $m$
is fixed in order to attain the best possible rates. This is
consistent with the procedures used in Section 2.2.

% We consider the best rate can be attained in Theorem \ref{tmn5}. The
% optimal $\bar{k}$ is
% \[
% \bar{k}\asymp T^{1/(2\lambda+1)}.
% \]
% Note that $\lambda\geq1$, the optimal $\bar{k}$ defined above is
% $o(T^{1/3})$ which matches the restriction of $\bar{k}$ in Theorem
% \ref{tmn5}. In the last section, we consider the case where the
% trend part is the linear combination of fixed observed factors and
% there is no correlation between the observed part and the latent
% part, which is the special case of subsection with $\lambda=\infty$.
% In the procedure of last section, we choose $\bar{k}$ as a fixed
% constant for the inference. Actually, the validity of such selection
% of $\bar{k}$ can be explained by Theorem \ref{tmn5}. In such a model
% setting, aim to obtain the best convergence rate of $\|\widehat{{\bf
% A}}-{\bf A}\|_2$, $\bar{k}\asymp 1$ via the above discussion of
% optimal $\bar{k}$.

Now we comment on the impact of $p$ on the convergence rate,
which depends critically on the factor strength
$\delta \in [0, 1]$ specified in Condition~\ref{as5}. To simplify
the notation, let $\kappa_1\asymp\kappa_2\asymp\kappa$ which is a
mild assumption in practice. Suppose
$p^{\delta}T^{(1-\lambda)/(2\lambda+1)}=o(1)$ and $\bar{k}\asymp
T^{1/(2\lambda+1)}$, Theorem \ref{tmn5} then reduces to
\[
\|\widehat{{\bf A}}-{\bf A}\|_2 = \left\{
                                  \begin{aligned}
                                    O_p(p^\delta
T^{-\lambda/(2\lambda+1)}), &~ \hbox{if $\kappa=o(p^{1-\delta})$;} \\
                                    O_p(p\kappa^{-1}T^{-\lambda/(2\lambda+1)}),
&~ \hbox{if
$p^{1-\delta}=o(\kappa)$.}\\
                                  \end{aligned}
                                \right.
\]
If $\kappa p^{-(1-\delta)}\rightarrow\infty$, there is an additional
factor $\kappa p^{-(1-\delta)}$ in the convergence rate of
$\|\widehat{{\bf A}}-{\bf A}\|_2$ than that under the setting
$\kappa p^{-(1-\delta)}\rightarrow0$, which implies that
$\|\widehat{{\bf A}}-{\bf A}\|_2$ converges to zero faster in the
case $\kappa = o(p^{1-\delta})$. The dimension $p$ must satisfy the
condition $p^{\delta}T^{(1-\lambda)/(2\lambda+1)}=o(1)$,
% If $\delta=0$, this condition
% is no restriction on $p$ to guarantee the consistency of
% $\widehat{{\bf A}}$.
which is automatically fulfilled when $\delta= 0$, i.e. the factors
are strong. However when the factors are weak in the sense
$\delta\neq 0$, $p$ can only be in the order
$p=o(T^{(\lambda-1)/\{(2\lambda+1)\delta\}})$ to ensure the
consistency in estimating the factor loading matrix.
% The following theorem shows the
% convergence rate of $p^{-1/2}\|\widehat{{\bf A}}\widehat{{\bf
% f}}_t-{\bf A}{\bf x}_t\|_2$.
\begin{tm}
\label{tmn6}  Let the condition of Theorem \ref{tmn5} hold.
% Under Conditions \ref{as1}-\ref{as2}, \ref{as4},
% \ref{as12}-\ref{as11}, suppose that $r$ is known, if
In addition, if $\|{\bf \Sigma}_\varepsilon\|_2$ is bounded as
$p\rightarrow\infty$, then
\[
p^{-1/2}\|\widehat{{\bf A}}\widehat{{\bf x}}_t-{\bf A}{\bf
x}_t\|_2=O_p(\|\widehat{{\bf A}}-{\bf
A}\|_2+p^{-1/2}+T^{-(2\lambda-1)/(4\lambda+2)}).
\]
\end{tm}

Comparing the above theorem with Theorem \ref{tm2}, it has one more
term $T^{(2\lambda-1)/(4\lambda+2)}$ in the convergence rate. When
the dimension $m$ is fixed and $\lambda=\infty$, it reduces to
Theorem \ref{tm2}.
% The only difference of the form between this theorem and Theorem
% \ref{tm2} is that the third term of the convergence rate determined
% by $\lambda$. In Theorem \ref{tm2} where the dimension of ${\bf
% z}_t$ is fixed and $\lambda=\infty$, the Theorem \ref{tmn6} with
% $\lambda=\infty$ is the same as Theorem \ref{tm2}.
On the other hand, we can also consider the model
(\ref{eq:nonlinear1}) with diverging number of regressors (i.e.,
$m\rightarrow\infty$). Noting Proposition \ref{pnn5} with
$\lambda=\infty$ and using the same argument of Theorem \ref{tm1},
it holds that
\[
\|\widehat{{\bf A}}-{\bf A}\|_2 = \left\{
                                  \begin{aligned}
                                    &~O_p\{\bar{k}^{-1}p^\delta(\bar{k}^{3/2}+m^{1/2})T^{-1/2}\}, \\
                                    &~~~~ \hbox{if $\kappa_2=o(p^{1-\delta})$ and $p^\delta(\bar{k}^{1/2}+m^{1/2})T^{-1/2}=o(1)$;} \\
                                    &~O_p\{\bar{k}^{-1}p\kappa_2\kappa_1^{-2}(\bar{k}^{3/2}+m^{1/2})T^{-1/2}\},\\
&~~~~ \hbox{if
$p^{1-\delta}=o(\kappa_1)$ and $p\kappa_2\kappa_1^{-2}(\bar{k}^{1/2}+m^{1/2})T^{-1/2}=o(1)$;}\\
                                  \end{aligned}
                                \right.
\]
provided that $m=o(T^{1/2})$ and $\bar{k}=o(T^{1/3})$. Theorem
\ref{tm1} can be regarded as the special case of this result with
fixed $\bar{k}$ and $m$. Note that the best convergence rate for
$\|\widehat{\bA}-\bA\|_2$ is attained under such setting if we
choose $\bar{k}\asymp m^{1/3}$.

\section{Numerical properties}

In this section, we illustrate the finite sample properties of the
proposed methods in two simulated models, one with a linear regression
term and one with a nonlinear regression term. For the linear model,
both stationary and nonstationary factors were employed. In each
model, we set the dimension of ${\bf y}_t$ at $p=100, 200, 400, 600,
800$ and the sample size $T=0.5p$, $p$, $1.5p$ respectively. For
each setting, 200 samples were generated.

\begin{table}[tbh]
\begin{center}
\captionstyle{center}
%\onelinecaptionsfalse
\caption{Relative frequency estimates of $P(\widehat{r}=r)$ for
Example 1 with stationary factors.}
%  with the true
% ${\bf D}$ (Oracle) or the estimator $\wh D$ (Real).}
\begin{tabular}{lll|ccccc}
  \hline
\hline
                  &          &        $p$       &      $100$         &    $200$         &     $400$        &  $600$  &   $800$     \\
   \hline
$\delta=0$ &  $\bD$  known   &      $T=0.5p$       & 0.700 &  0.960
& 0.990 & 0.995 & 1

 \\
           &             &        $T=p$            &  0.900 &  0.985  & 1 & 1& 1\\
           &             &        $T=1.5p$         &  0.980 & 1& 1 & 1 & 1 \\
\cline{2-8}
           &  $\bD$ unknown      &        $T=0.5p$ & 0.615   & 0.940 & 0.990 & 0.995 & 1 \\
           &             &         $T=p$           &  0.865  & 0.985 & 1  & 1 & 1 \\
           &             &         $T=1.5p$        &  0.960  & 1 & 1 & 1& 1 \\
\hline
$\delta=0.5$ &   $\bD$ known   &        $T=0.5p$   &  0.105 & 0.805 & 0.950 & 0.805 & 0.930\\
           &             &        $T=p$            &  0.285  & 0.880 & 0.940 & 1 & 0.975\\
           &             &        $T=1.5p$         &   0.895 &   0.975& 1& 1& 1 \\
\cline{2-8}
           &  $\bD$ unknown      &        $T=0.5p$ &  0.065 & 0.765  & 0.930 & 0.780& 0.910 \\
           &             &         $T=p$           &  0.280 & 0.880  &  0.940  & 1 & 0.975\\
           &             &         $T=1.5p$        &   0.870 & 0.975  & 0.995  & 1 & 1 \\
 \hline
\hline
\end{tabular}
 \label{table1}
\end{center}
\end{table}

\begin{figure}[tbh]
\subfigure{\includegraphics[scale=0.4]{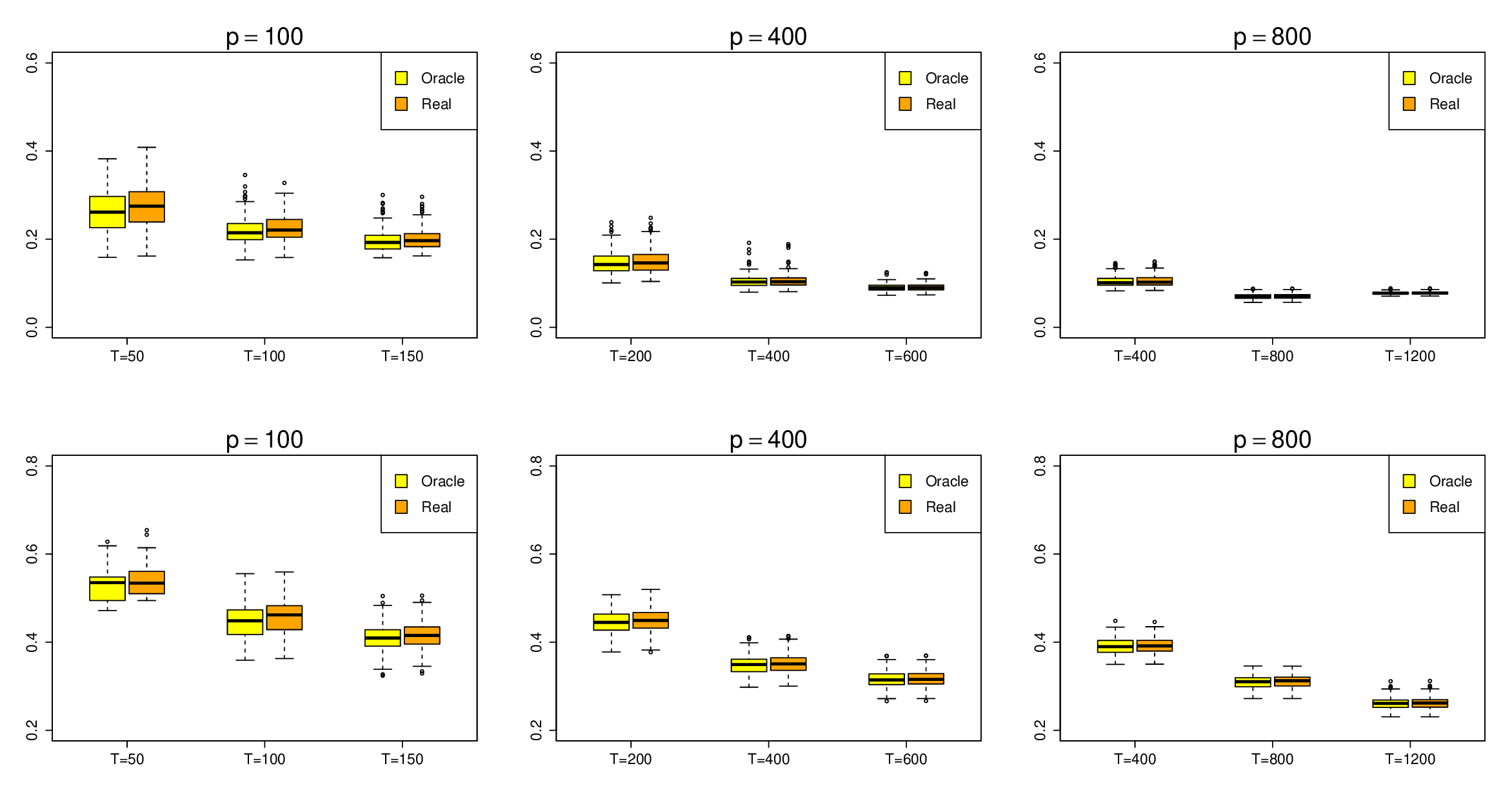}}\caption{Boxplots
of $\{D(\mathcal{M}(\widehat{{\bf A}}),\mathcal{M}({\bf A}))\}$ for
Example 1 with stationary factor, and $\delta=0$ (3 top panels) and
$\delta=0.5$ (3 bottom panels). Errors obtained using true $\bD$ are
marked with `oracle', and using $\wh \bD$ are marked with `real'.
\label{fig4}}
\end{figure}

\noindent {\bf Example 1}. Consider the linear model ${\bf y}_t=
{\bf D}{\bf z}_t + {\bf A}{\bf x}_t +\bvarepsilon_t$, in which ${\bf
z}_t$ follows the VAR(1) model:
\begin{equation} \label{ex1}
{\bf z}_t=\left(
            \begin{array}{cc}
             5/8  & {1}/{8} \\
               {1}/{8}  &  {5}/{8} \\
            \end{array}
          \right){\bf z}_{t-1}+{\bf e}_{t},
\end{equation}
where ${\bf e}_{t}\sim N({\bf 0},{\bf I}_2)$. Let ${\bf D}$ be a
$p\times 2$ matrix of which the elements were generated independently from
the uniform distribution $U(-2,2)$,  ${\bf x}_t$ be $3\times 1$
VAR(1) process with independent $N({\bf 0},{\bf I}_3)$ innovations
and the diagonal autoregressive coefficient matrix with 0.6, -0.5
and 0.3 as the main diagonal elements. This is a stationary factor
process with $r=3$ factors.  The elements of ${\bf A}$ were drawn
independently from $U(-2,2)$
 resulting a strong factor case with
$\delta=0$. Also we considered a weak factor case with $\delta=0.5$
for which randomly selected $p-\lfloor{ p^{1/2} \rfloor}$ elements
in each column of ${\bf A}$ were set to 0.
% Then, ${\bf A}$ was transformed such that ${\bf A}^{\T}{\bf A}={\bf I}_r$.
Let $\bvarepsilon_t$ be independent and $N({\bf 0},{\bf I}_p)$. To
show the impact of the estimated coefficients matrix ${\wh \bD}$ on
the estimation for the factors, we also report the results from
using the true ${\bf D}$. We report the results with $\bar{k}=1$
only, since the results with $1 \le \bar{k}\le 10$ are similar. The
relative frequency estimates of $P(\widehat{r}=r)$ are reported in
Table \ref{table1}. It shows that the defect in estimating $r$ due
to the errors in estimating $\bD$ is almost negligible.
Fig.\ref{fig4} displays the boxplots of the estimation errors
$\{D(\mathcal{M}(\widehat{{\bf A}}),\mathcal{M}({\bf A}))\}$. Again
the performance with the estimated coefficient matrix $\wh \bD$ is
only slightly worse than that with the true $\bD$. When the factors
are weaker (i.e. when $\delta =0.5$), it is harder to estimate both
the number of factors and the factor loading space. All those
findings are in line with the asymptotic results presented in
Section 2.3.

\begin{figure}[tbh]
\subfigure{\includegraphics[scale=0.4]{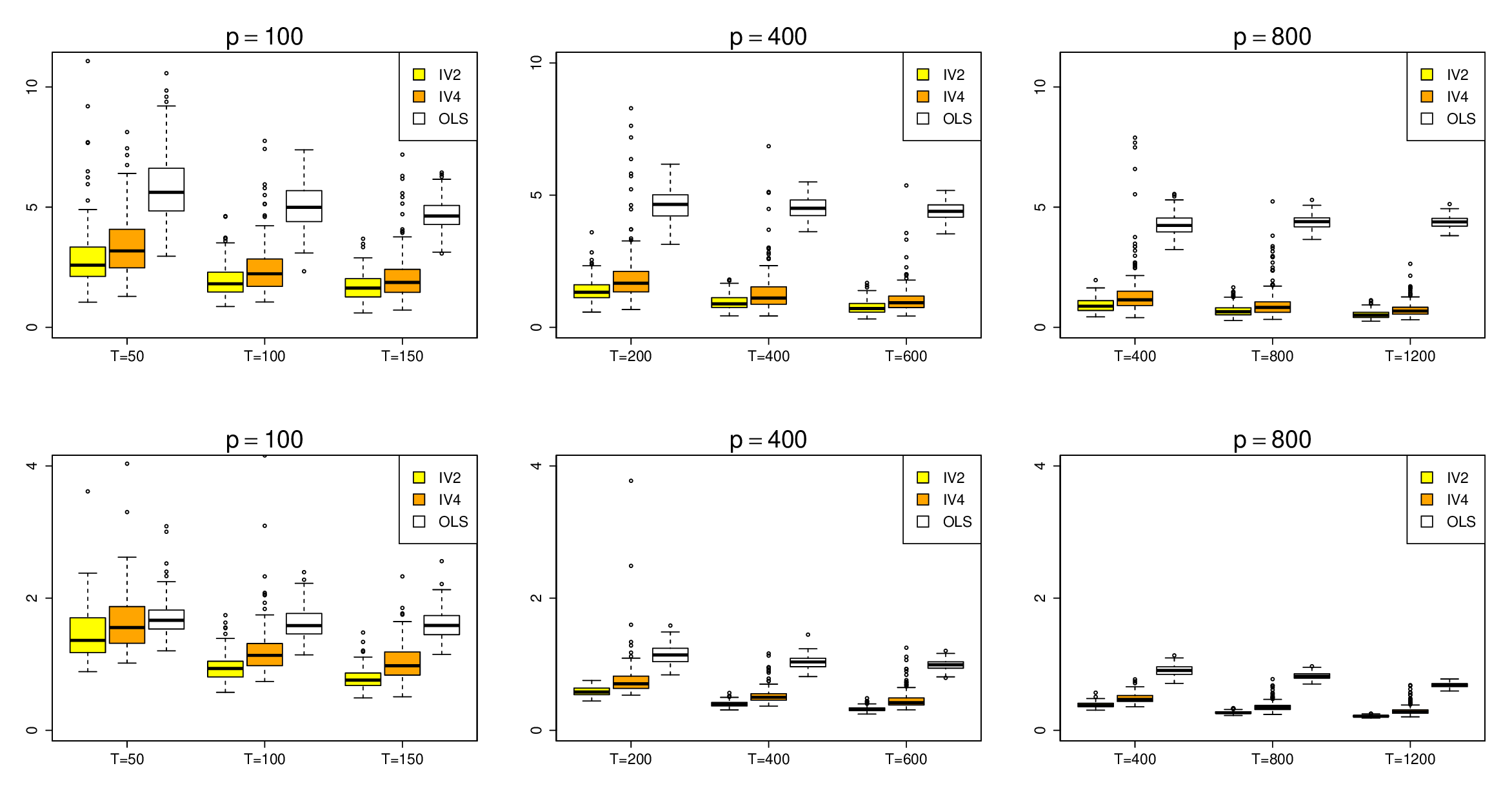}}
\caption{Boxplots of  $p^{-1/2}\|\widehat{\bf D}-{\bf D}\|_F$ for
Example 1 with endogeneity,  and $\delta=0$ (3 top panels) and
$\delta=0.5$ (3 bottom panels). \label{fig5}}

\end{figure}
\begin{table}[tbh]
\begin{center}

 \captionstyle{center}
\onelinecaptionsfalse\caption{Relative frequency estimates of
$P(\widehat{r}=r)$ for Example 1 with endogeneity.}
% using instrument variable (IV) or ordinary least squares (OLS) methods.}
\begin{tabular}{lll|ccccc}
  \hline
  \hline
&          &        $p$       &      $100$         &    $200$         &     $400$        &  $600$  &   $800$     \\
   \hline
 &    IV$2$             &        $T=0.5p$          &   0.660  &  0.885  &  0.995  &  0.995  & 1 \\
           &             &        $T=p$            &  0.855   &  0.990  & 1       & 1       &1 \\
           &             &        $T=1.5p$         &  0.960   &  1      & 1       &  1       & 1 \\
\cline{2-8}
           &    IV$4$       &        $T=0.5p$      &  0.590   &   0.865 & 0.975   & 0.970   & 0.965 \\
$\delta=0$           &             &        $T=p$  & 0.845    &   0.970 & 0.970   & 0.990   & 0.985\\
           &             &        $T=1.5p$         & 0.930    &   0.990 &  0.980  & 0.970        & 0.975  \\
\cline{2-8}
           &    OLS      &        $T=0.5p$         &  0.580   &  0.865  & 0.990   & 1       & 1\\
           &             &         $T=p$           &  0.855   &  0.980  & 1       & 1       & 1\\
           &             &         $T=1.5p$        &  0.945   &   1     & 1       &  0.995       & 1\\
\hline
 &   IV$2$      &        $T=0.5p$                  &   0.280  &  0.665  & 0.625   & 0.630   & 0.620 \\
           &             &        $T=p$            &  0.600   &  0.715  &
           0.980 & 1 & 1
\\
           &             &        $T=1.5p$         &  0.550   &   0.980 & 0.990   & 1 & 1 \\
\cline{2-8}
           &    IV$4$       &        $T=0.5p$      & 0.205    &   0.570 & 0.550  & 0.605  & 0.550
\\
$\delta=0.5$       &    &        $T=p$             & 0.510 & 0.650     & 0.890  & 0.960  & 0.925\\
           &             &        $T=1.5p$         & 0.580 &  0.915    & 0.950  & 0.935 & 0.940 \\
\cline{2-8}
           &     OLS     &        $T=0.5p$         & 0.225 &   0.405  & 0.635   & 0.625 & 0.640\\
           &             &         $T=p$           & 0.535 & 0.705    & 1       &  1    & 0.995\\
           &             &         $T=1.5p$        & 0.630 & 0.955    & 1      & 1 & 0.995\\
\hline \hline

\end{tabular}
\label{table2}
\end{center}
\end{table}
\begin{figure}[tbh]
\subfigure{\includegraphics[scale=0.4]{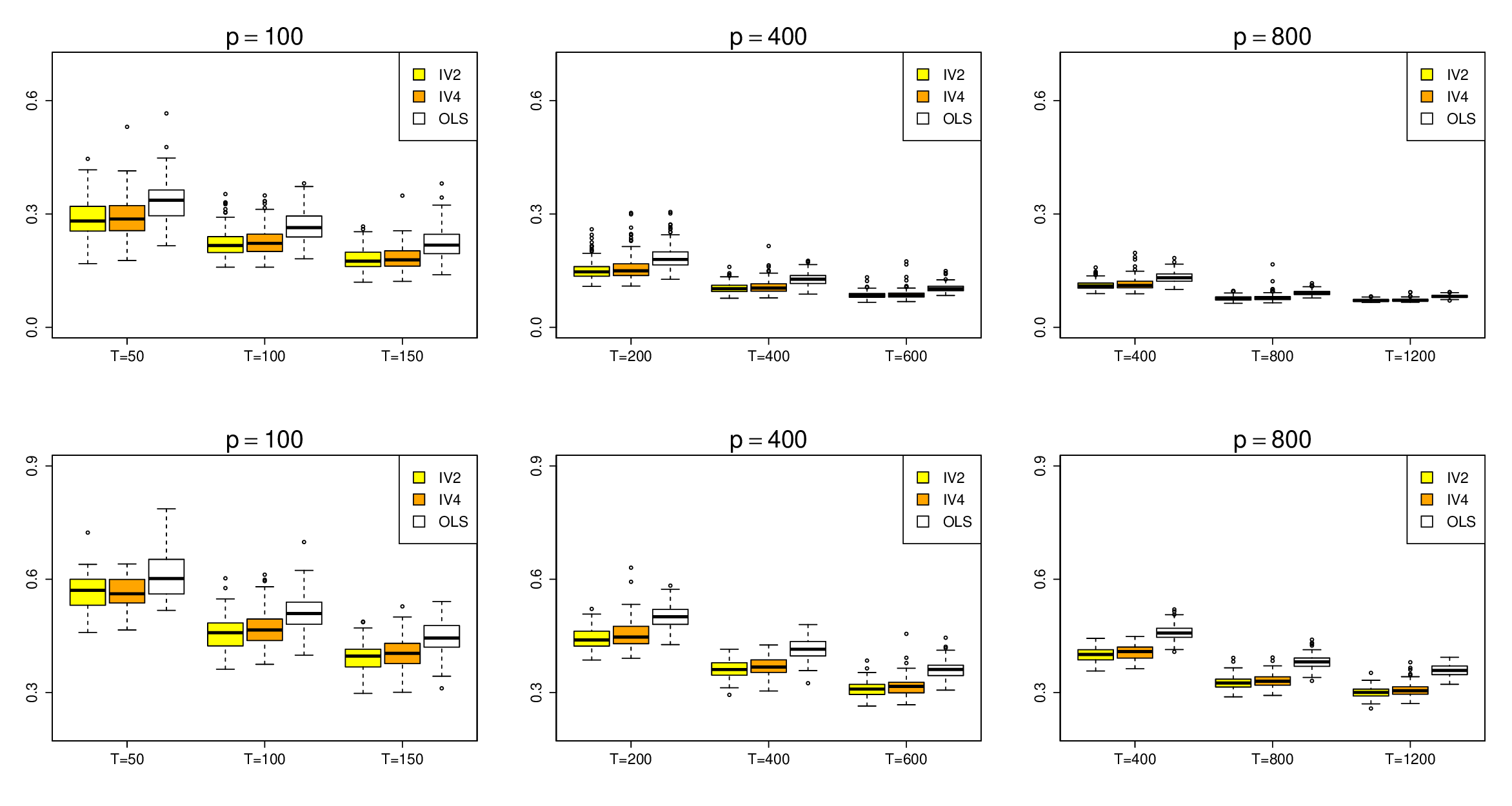}}\setlength{\abovecaptionskip}{0pt}
\setlength{\belowcaptionskip}{10pt}
% \captionstyle{center} \onelinecaptionsfalse
\caption{Boxplots of  $\{D(\mathcal{M}(\widehat{{\bf
A}}),\mathcal{M}({\bf A}))\}$ for Example 1 with endogeneity, and
$\delta=0$ (3 top panels) and $\delta=0.5$ (3 bottom panels).
\label{fig6}}

\end{figure}

Now we consider the case with the endogeneity. To this end, we
changed the definition for the regressor process $\bz_t$ in the
above setting. Instead of (\ref{ex1}), we let
\[
z_{1,t}=0.1x_{1,t}+0.1u_t+0.1u_t^2, \quad
z_{2,t}=0.1x_{2,t}-0.1u_t+0.1u_t^2,
\]
where $ u_t $ is an AR(1) process defined by $
u_t=0.5u_{t-1}+\epsilon_t$ and $\epsilon_t\sim N(0,1)$. The ordinary
least squares estimator of ${\bf D}$ is no longer consistent now. We
employ  two different instrument variables ${\bf
w}_t=(u_t,u_t^2)^{\T}$ and $\widetilde{{\bf
w}}_t=(u_t,u_t^2,u_t^3,u_t^4)^{\T}$, as they are correlated with
${\bf z}_t$ but uncorrelated with ${\bf x}_t$ and $\bvarepsilon_t$.
The estimation error for $\widehat{{\bf D}}$ is measured by the
normalized Frobenius norm $p^{-{1}/{2}}\|\widehat{{\bf D}}-{\bf
D}\|_{F}$.  Setting ${\bf R}={\bf I}_2$ for ${\bf w}_t$ and the
elements of $\wt {\bf R}$ are generated from $U(-2,2)$ for
$\widetilde{{\bf w}}_t$ in (\ref{eq:ive}), we computed first both
the ordinary least squares (OLS) estimates and the instrument
variable method (IV) estimates for $\bD$, and then the estimates for
the number of factors $r$ and the factor loading matrix $\bA$ based
on, respectively, the two sets of residuals resulted from the two
regression estimation methods. The results are reported in
Figs.\ref{fig5} and \ref{fig6} and Table \ref{table2}   where
IV2 and IV4 represent the estimation using ${{\bf w}}_t$ and
$\widetilde{{\bf w}}_t$ respectively. Those simulation results
reinforce the findings in Section 3, which indicate that the
existence of the endogeneity has no impact in identifying and in
estimating the factor loading space. More precisely, Fig.\ref{fig5}
shows that the errors $p^{-{1}/{2}}\|\widehat{{\bf D}}-{\bf
D}\|_{F}$ for the OLS method are unusually large, as it effectively
estimates $\bD^\star$ in (\ref{d1}) instead of $\bD$. On the other
hand, the IV method provides accurate estimates for $\bD$. However
the differences of the two methods on the subsequent estimation for
the number of factors $r$ and the factor loading space $\calM(\bA)$
are small; see Table \ref{table2} and Fig.\ref{fig6}. Since the IV
method uses extra information, %${\bf w}_t$,
 it tends to offer
slightly better performance. Nevertheless Table \ref{table2}
indicates  that this improvement in estimating $r$ is almost
negligible. Also, the results are not sensitive to the choice
of $\bR$ as long as the instrument variables are properly selected.

\begin{table}[htb]
\begin{center}
% \captionstyle{center} \onelinecaptionsfalse
\caption{Relative frequency estimates of $P(\widehat{r}=r)$ for
Example 1 with nonstationary factors.}
% linear model as if we have known ${\bf D}$ (Oracle) or not (Real).}

\vspace{1ex}

\begin{tabular}{lll|ccccc}
  \hline
  \hline
                  &          &        $p$       &      $100$         &    $200$         &     $400$        &  $600$  &   $800$     \\
   \hline
$\delta=0$ &  $\bD$  known   &        $T=0.5p$     &  0.155 & 0.525 & 0.855 & 0.925 & 0.970\\
           &             &        $T=p$            &  0.465 & 0.800 & 0.940 & 0.990 & 0.990\\
           &             &        $T=1.5p$         &  0.625 & 0.890 & 0.995 & 0.985 & 1\\
\cline{2-8}
           &  $\bD$ unknown      &        $T=0.5p$ &  0.110 & 0.525 &0.835 &  0.920 & 0.970\\
           &             &         $T=p$           &  0.430 & 0.780 &0.940 & 0.990  & 0.990\\
           &             &         $T=1.5p$        &  0.595 & 0.890 &0.995 & 0.985 & 1\\
\hline
$\delta=0.5$ &   $\bD$ known   &        $T=0.5p$   &  0    &  0.070 & 0.175 & 0.385 & 0.525  \\
           &             &        $T=p$            & 0.025 & 0.235  & 0.535 & 0.705 & 0.765 \\
           &             &        $T=1.5p$         &  0.145& 0.475  & 0.740 & 0.815 & 0.860 \\
\cline{2-8}
           &  $\bD$ unknown      &        $T=0.5p$ &   0   & 0.055 & 0.160 & 0.380 & 0.520\\
           &             &         $T=p$           &  0.025& 0.215 & 0.520 & 0.685 & 0.760\\
           &             &         $T=1.5p$        &  0.125& 0.465 & 0.740 & 0.805 & 0.850\\
\hline \hline

\end{tabular}
\label{table3}
\end{center}
\end{table}

\begin{figure}[htb]
\subfigure{\includegraphics[scale=0.4]{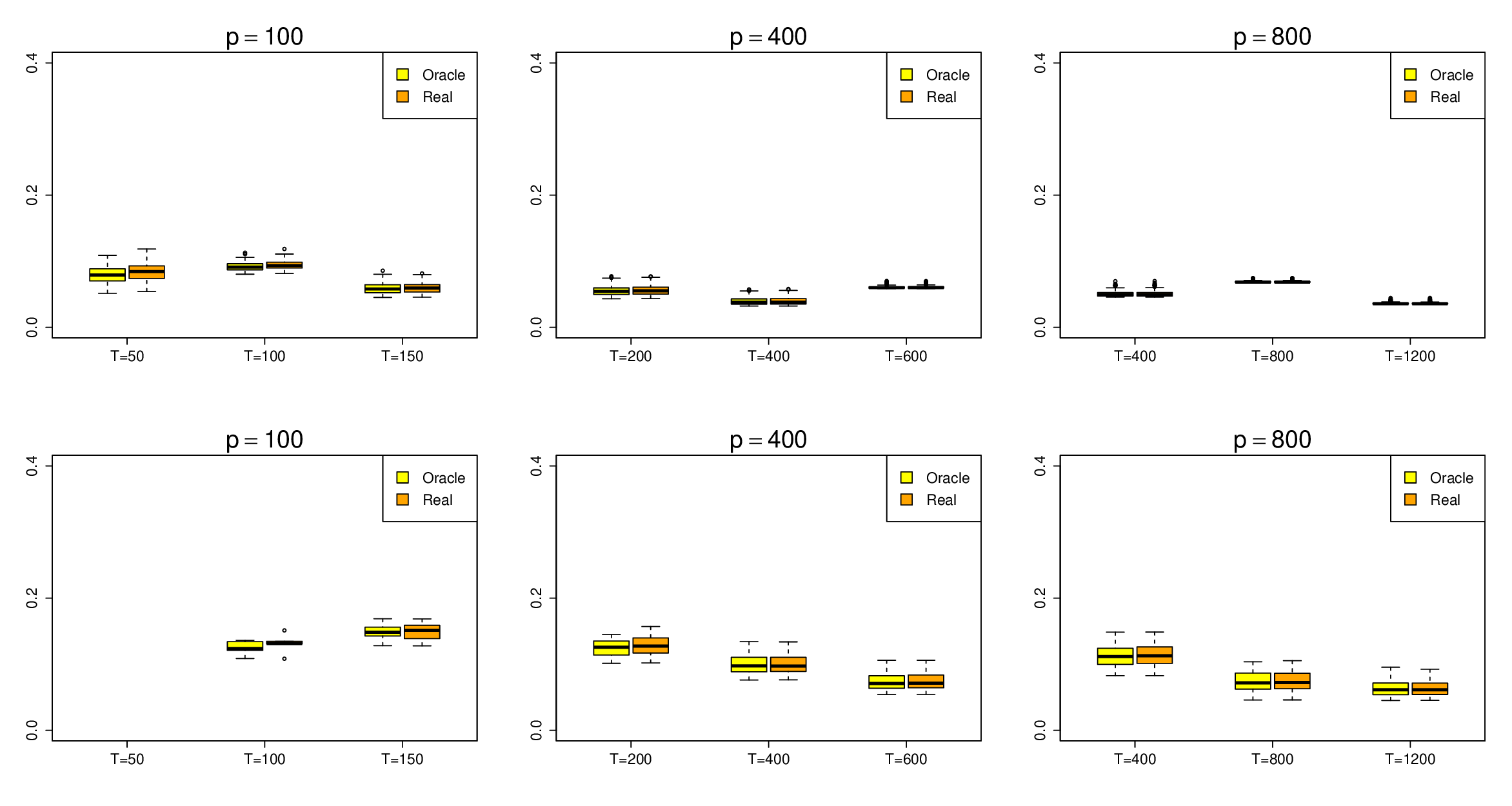}}
 \setlength{\abovecaptionskip}{0pt}
\setlength{\belowcaptionskip}{10pt}
% \captionstyle{center} \onelinecaptionsfalse
\caption{Boxplots of $\{D(\mathcal{M}(\widehat{{\bf
A}}),\mathcal{M}({\bf A}))\}$ for Example 1  with nonstationary
factors,  and $\delta=0$ (3 top panels) and $\delta=0.5$ (3 bottom
panels). Errors obtained using true $\bD$ are marked with `oracle',
and using $\wh \bD$ are marked with `real'. \label{fig7}}
\end{figure}

Now we consider the model with nonstationary factors
$\bx_t=3(x_{1,t}, x_{2,t}, x_{3,t})^{\T}$:
\begin{equation} \label{ex11}
x_{1,t}-2t/T=0.8(x_{1,T-1}-2t/T)+e_{1,t}, \quad x_{2,t}=3t/T, \quad
x_{3,t}=x_{3,t-1}+\sqrt{\frac{10}{T}}e_{3,t},
%~~\text{with}~~X_{0,3}\sim N(0,1),
\end{equation}
where $e_{j,t}$ are independent and $N(0, 1)$. The other settings
are the same as the first part of this example. The results are
reported in Table \ref{table3} and Fig.\ref{fig7}. The patterns are
similar to those in Table \ref{table1} and Fig.\ref{fig4}, except
that for a fixed $p$, the performance does not necessarily improve
when the sample size $T$ increases; see Fig.\ref{fig7}. This is due
to the nonstationary nature of the factors defined in (\ref{ex11}):
new observations bring in the information on the new and
time-varying underlying structure as far as the factor processes are
concerned.
%We also notice that now the accuracy of the estimation is not as good as
%that in the stationary cases.

\begin{table}[tbh]
\begin{center}
% \captionstyle{center}
%\onelinecaptionsfalse
\caption{Relative frequency estimates of $P(\widehat{r}=r)$ for
Example 2 (with nonlinear regression).  }
% as if we have known ${\bf g}(\cdot)$ (Oracle) or not (Real).}
\begin{tabular}{lll|ccccc}
  \hline
  \hline
                  &          &        $p$       &      $100$         &    $200$         &     $400$        &  $600$  &   $800$     \\
   \hline
$\delta=0$ &    ${\bf g}$ known   &       $T=0.5p$ & 0.780  & 0.865 & 0.965 & 0.975 & 0.985\\
           &             &        $T=p$            & 0.840  & 0.920 & 0.990 & 1     & 1\\
           &             &        $T=1.5p$         & 0.820  & 0.990 &1      &1      & 1\\
\cline{2-8}
           &   ${\bf g}$ unknown      &   $T=0.5p$ & 0.750  & 0.860 & 0.955 & 0.975 & 0.980\\
           &             &         $T=p$           & 0.830  & 0.890 & 0.990 & 1     & 1\\
           &             &         $T=1.5p$        & 0.780  & 0.990 & 1     &1      & 1\\
\hline
$\delta=0.5$ &    ${\bf g}$ known   & $T=0.5p$     & 0.270 & 0.665 & 0.725 & 0.430 & 0.650\\
           &             &        $T=p$            & 0.390 & 0.700  & 0.850  & 0.810 & 0.800 \\
           &             &        $T=1.5p$         & 0.390 & 0.720 & 0.885 & 0.960 &1\\
\cline{2-8}
           &   ${\bf g}$ unknown      &$T=0.5p$    & 0.260 & 0.625 & 0.665 & 0.390 & 0.600 \\
           &             &         $T=p$           & 0.390 & 0.655 & 0.760 & 0.810 & 0.795\\
           &             &         $T=1.5p$        & 0.335 & 0.700 & 0.875 & 0.950 &1\\
\hline \hline

\end{tabular}
\label{table4}
\end{center}
\end{table}

\begin{figure}[htb]
\subfigure{\includegraphics[scale=0.4]{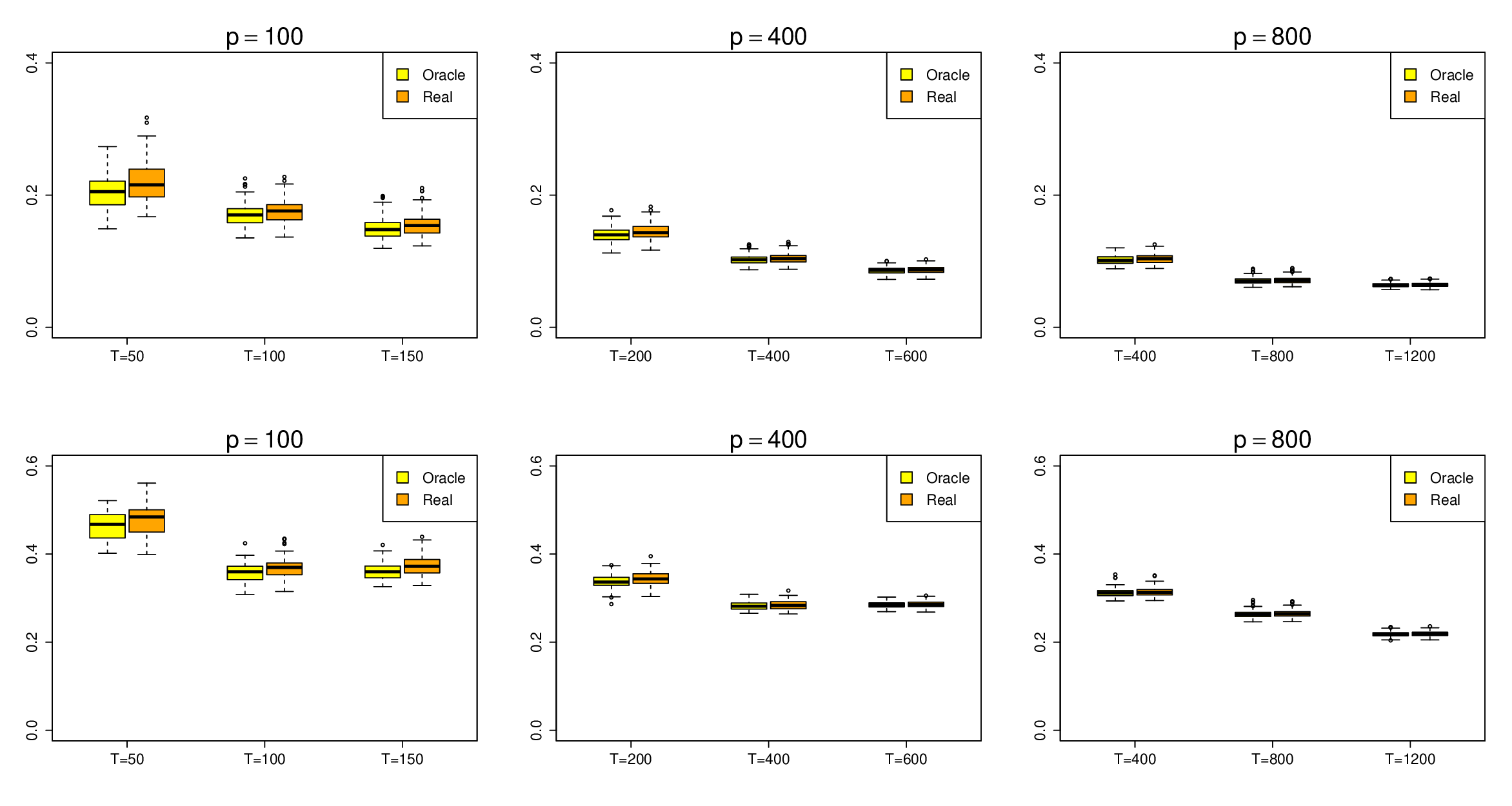}}
\caption{Boxplots of $\{D(\mathcal{M}(\widehat{{\bf
A}}),\mathcal{M}({\bf A}))\}$ for Example 2  with nonlinear
regression,  and $\delta=0$ (3 top panels) and $\delta=0.5$ (3
bottom panels). Errors obtained using true ${\bf g}$ are marked with
`oracle', and using $\wh {\bf g}$ are marked with `real'.
\label{fig8}}
\end{figure}

\noindent {\bf Example 2}. We now consider a model with nonlinear
regression function. Let $\bu_t = u_t$ be a univariate AR(1) process
defined by $u_t=0.5u_{t-1}+e_t$ with independent $N(0,1)$
innovations $e_t$. The nonlinear regression function ${\bf
 g}(u_t)=(g_1(u_t),\ldots,g_p(u_t))^{\T}$ was defined as
 \[
 g_i(u_t)=\frac{\exp(\alpha_i^{(1)}u_t)}{1+\exp(\alpha_i^{(1)}u_t)},~~i=1,\ldots,\frac{p}{2}~~\text{and}~~g_i(u_t)=\sin(\alpha_i^{(2)}u_t),~~i=\frac{p}{2}+1,\ldots,
 p,
 \]
where the parameters $\alpha_i^{(1)}$ were drawn independently from
$N(0,4)$,  and $\alpha_i^{(2)}$ were drawn independently from
 $U(-2,2)$ respectively. We used the same $\bA, \bx_t $ and $\bve_t$ as
in the first part of Example 1.

We used the polynomial expansion to approximate ${\bf
 g}(u_t)$, i.e.
$g_i(u_t)\approx \sum_{j=1}^{m}d_{i,j}l_j(u_t)$ with $l_j(
u_t)=u_t^{j-1}$, where the order $m$ was set  as $\lfloor
2T^{1/5}\rfloor$. We obtained $\hat{d}_{i,j}$ by the least square
estimation. Put $\widehat{{\bf
g}}(u_t)=(\hat{g}_1(u_t),\ldots,\hat{g}_p(u_t))^{\T}$ for
$\hat{g}_i(u_t)= \sum_{j=1}^{m}\hat{d}_{i,j}l_j(u_t)$. The residuals
$\widehat{\beeta}_t={\bf y}_t-\widehat{{\bf g}}({\bf u}_t)$ were
then used to estimate the latent factors. We set $\bar{k}=\lfloor
2T^{1/5}\rfloor$; see Theorem \ref{tmn5}. The simulation results are
reported in Table \ref{table4} and Fig.\ref{fig8}, which present
similar patterns as in the first part of Example 1.

\section{Real data analysis}
We illustrate our method by modeling the daily returns of $123$
stocks from $2$ January 2002 to 11 July 2008. The stocks were
selected among those contained in the S\&P$500$ which were traded
everyday during this period. The returns were calculated based on
the daily close prices. We have in total $T=1642$ observations with
the dimension $p=123$. This data has been analyzed in
\cite{LamYao_AOS_2012}. They identified two factors under a pure
factor model setting, i.e. model (\ref{eq:nonlinear1}) with ${\bf
z}_t\equiv \bf0$. Furthermore the estimated factor loading space
contains the return of the S\&P$500$. Hence it can be regarded as
one of the two factors. Since the S\&P$500$ index is often viewed as
a proxy of the market index,
% Usually, the daily return of each stock is closely related to the S\&P$500$ index.
it is reasonable to take its return  as a known factor $z_t$ in our
model (\ref{eq:nonlinear1}).
% We standardized the data firstly.
We calculated the  ordinary least square estimator for the
regression coefficient matrix ${\bf D}$ which is now a $123\times 1$
vector with each element representing the impact of the S\&P$500$
index to the return of the corresponding stock. As all the estimated
elements are positive, indicating the positive correlations between
the returns of market index and the those 123 stocks.

% We then identify the number of
% latent factors. The parameter $\bar{k}$ was chosen from $1$ to $20$
% respectively to determine $\widehat{\bf M}$ defined in
% (\ref{eq:hatm}). The estimations for the number of the latent
% factors based on different selections of $\bar{k}$ were the same.
% Without loss of generality, we set $\bar{k}=1$ in our following
% analysis.
\begin{figure}
\centering \subfigure[The eigenvalues of $\widehat{{\bf M}}$
]{\includegraphics[scale=0.4]{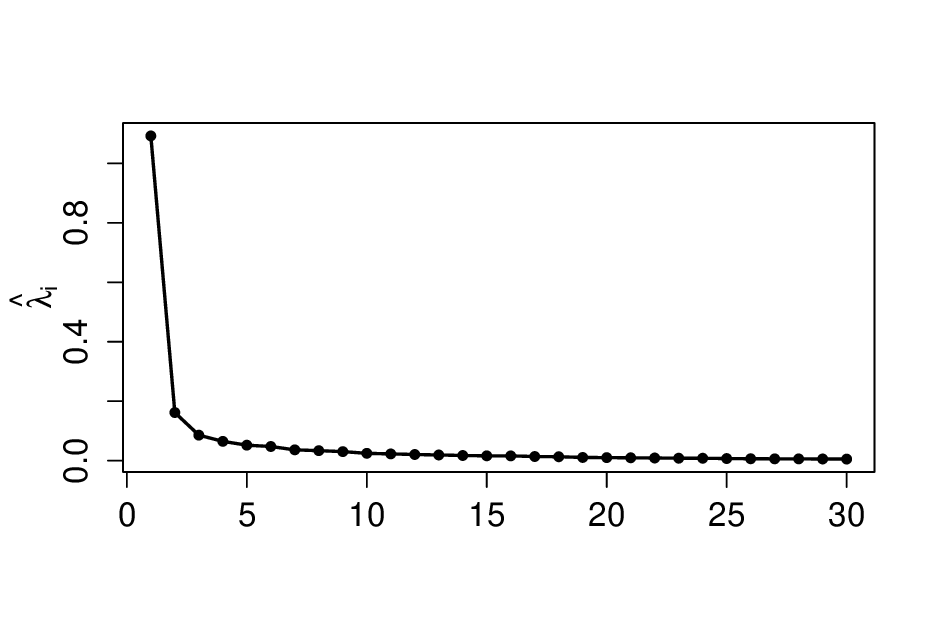}}
\subfigure[Ratio of the eigenvalues of $\widehat{{\bf M}}$ ]
{\includegraphics[scale=0.4]{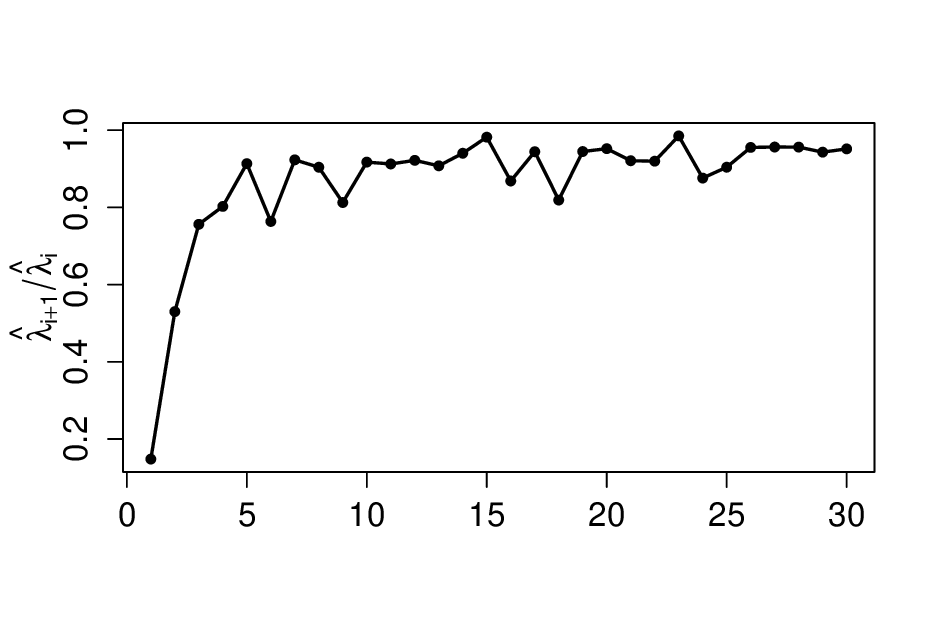}}
\subfigure[The estimated latent factor ]
{\includegraphics[width=0.8\linewidth]{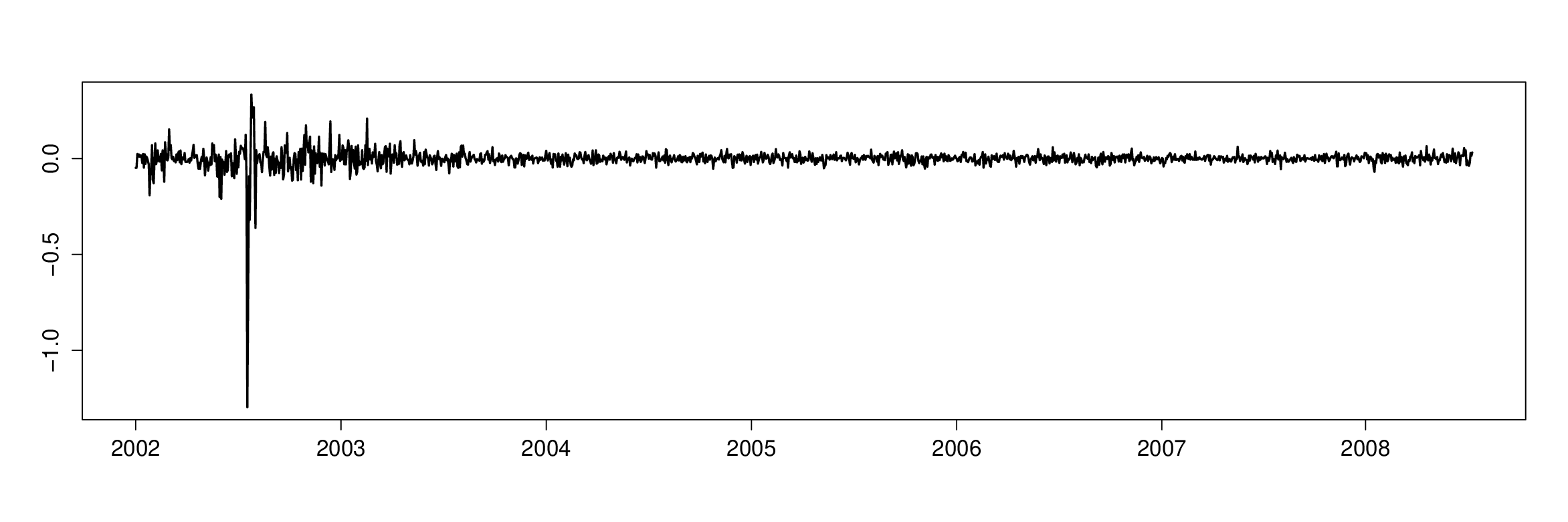}}
\subfigure[S\&P $500$ returns]
{\includegraphics[width=0.8\linewidth]{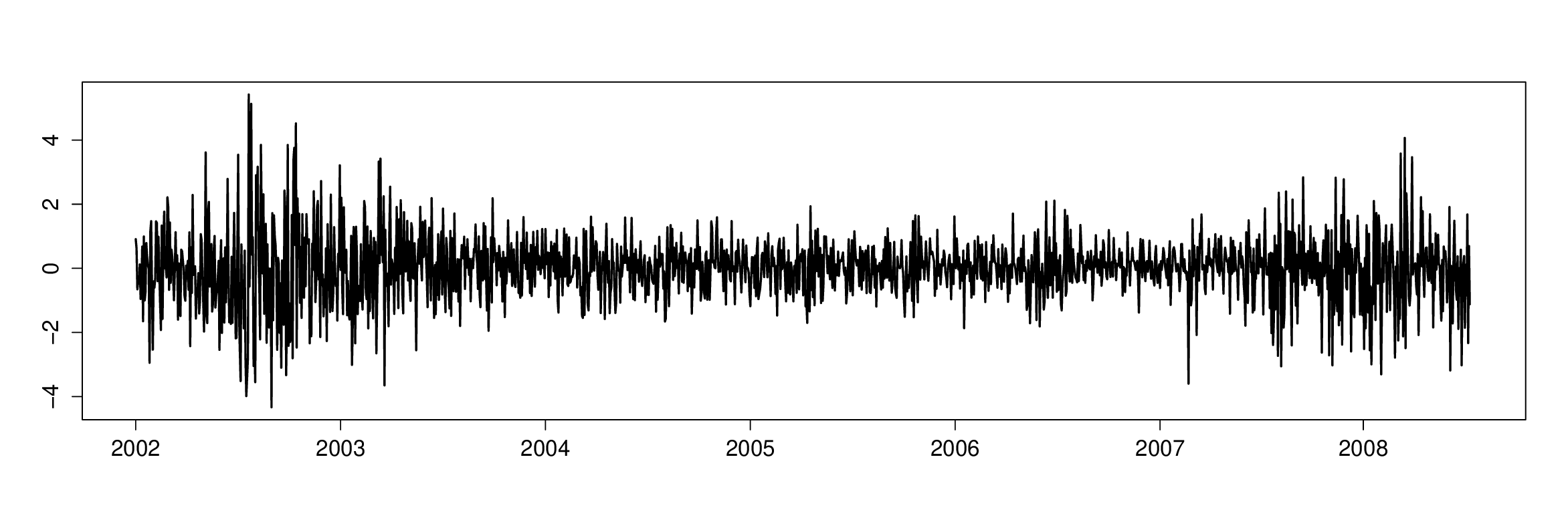}}
 \label{figure1}
\caption{The estimated eigenvalues (multiplied by $10^6$), the ratio
of eigenvalues, the estimated latent factor and the S\&P$500$ returns in
2 January 2002 -- 11 July 2008.}

\end{figure}

\begin{figure}
\centering
\subfigure{\includegraphics[width=.27\linewidth]{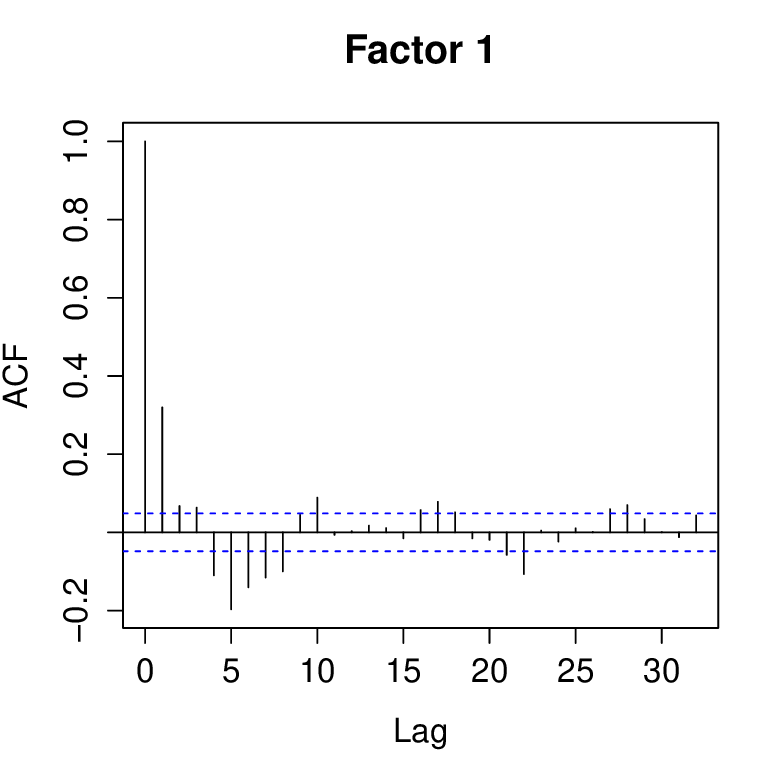}}
\subfigure{\includegraphics[width=.27\linewidth]{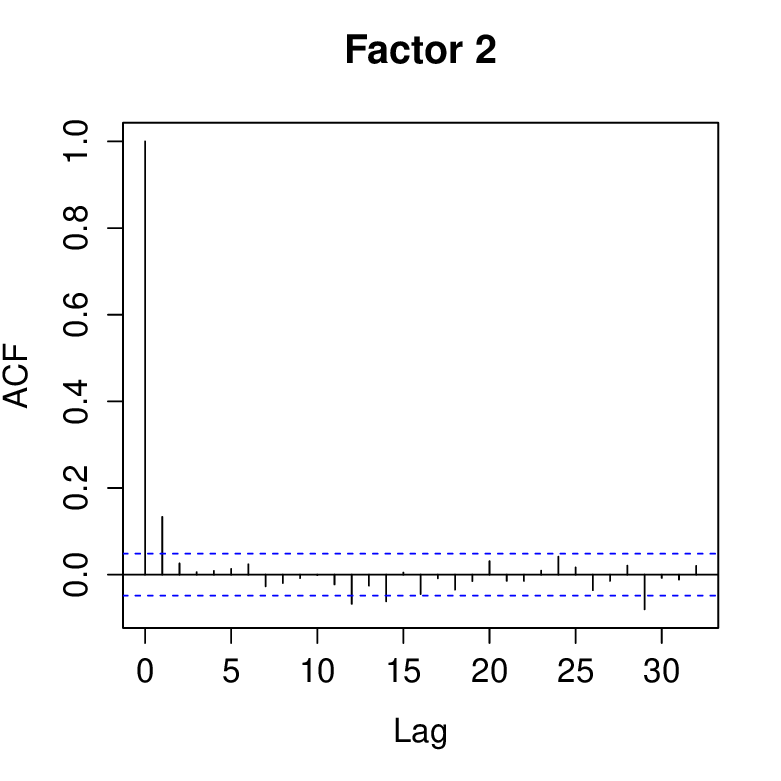}}
\subfigure {\includegraphics[width=.27\linewidth]{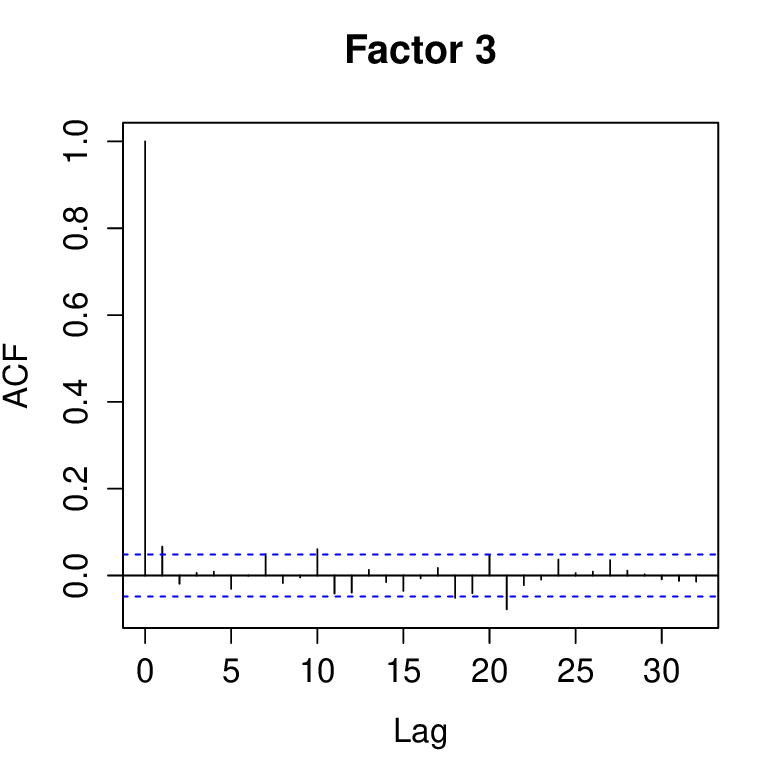}}
%\subfigure{\includegraphics[width=.3\linewidth]{figure/pacf1.eps}}
%\subfigure{\includegraphics[width=.3\linewidth]{figure/pacf2.eps}}
%\subfigure {\includegraphics[width=.3\linewidth]{figure/pacf3.eps}}
\caption{The ACFs of the first three estimated factors.}
\label{figure2}
\end{figure}

Fig.\ref{figure1}  displays the first $30$ eigenvalues of
$\widehat{{\bf M}}$, defined as in (\ref{eq:hatm}) with $\bar k=1$,
 sorted in the descending order. The ratio of
$\widehat{\lambda}_{i+1}/\widehat{\lambda}_{i}$ in the right panel
indicates that there is only one latent factor. Varying $\bar k$
between 1 to 20 did not alter this result. Fig.\ref{figure1}(c)
shows that the sparks of the estimated factor process occur around
$22$ July, 2002, which is consistent with the oscillations of
S\&P500 index, although the S\&P$500$ are less volatile.
%since the index is an average of $500$ stocks that implies it should
%be more stable.
The autocorrelations of the estimated factors
$\widehat{\gamma}_j^{\T}({\bf y}_t-\widehat{\bf D}{\bf z}_t)$, where
$\widehat{\gamma}_j$ is the unit eigenvector of $\widehat{{\bf M}}$
corresponding to its $j$th largest eigenvalue, are plotted in
Fig.\ref{figure2} for $j=1, 2, 3$. The autocorrelations of the first
factor is significant non-zero. On the other hand, there are hardly
any significant non-zero autocorrelations for both the second and
the third factors.
% indicates that the first factor contains almost all of the information,
% the other factors are not significant against white noise. Hence,
% setting $\widehat{r}=1$ is reasonable.

\begin{figure}
\centering
\subfigure{\includegraphics[scale=0.5]{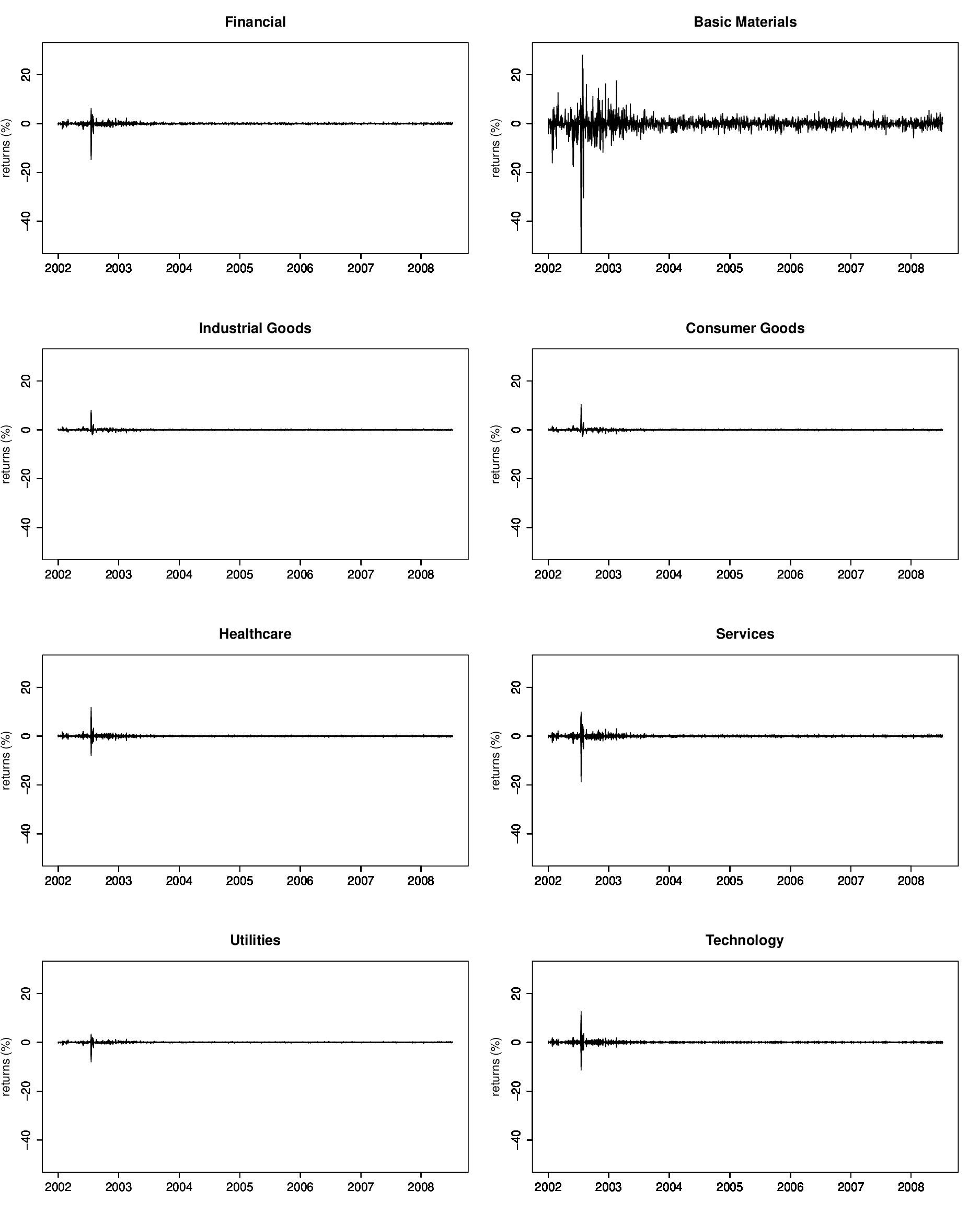}}   %1200*1500
\caption{The estimated latent part ${\bf A}{\bf x}_t$ across
different sectors.} \label{figure3}
\end{figure}

\begin{figure}
\centering \subfigure[The eigenvalues of $\widehat{{\bf M}}$
]{\includegraphics[scale=0.4]{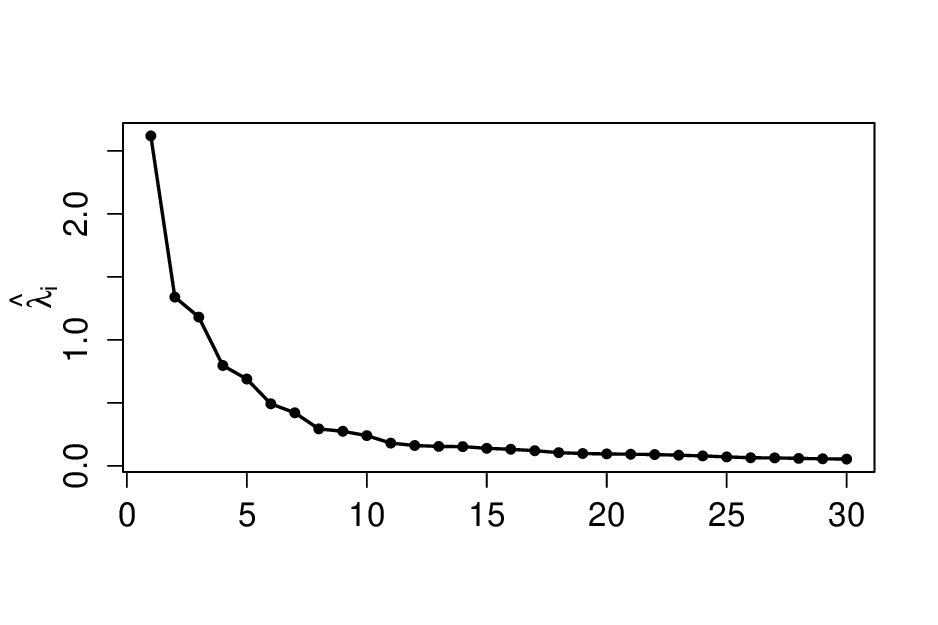}}        % 600*400
\subfigure[Ratio of the eigenvalues of $\widehat{{\bf M}}$ ]
{\includegraphics[scale=0.4]{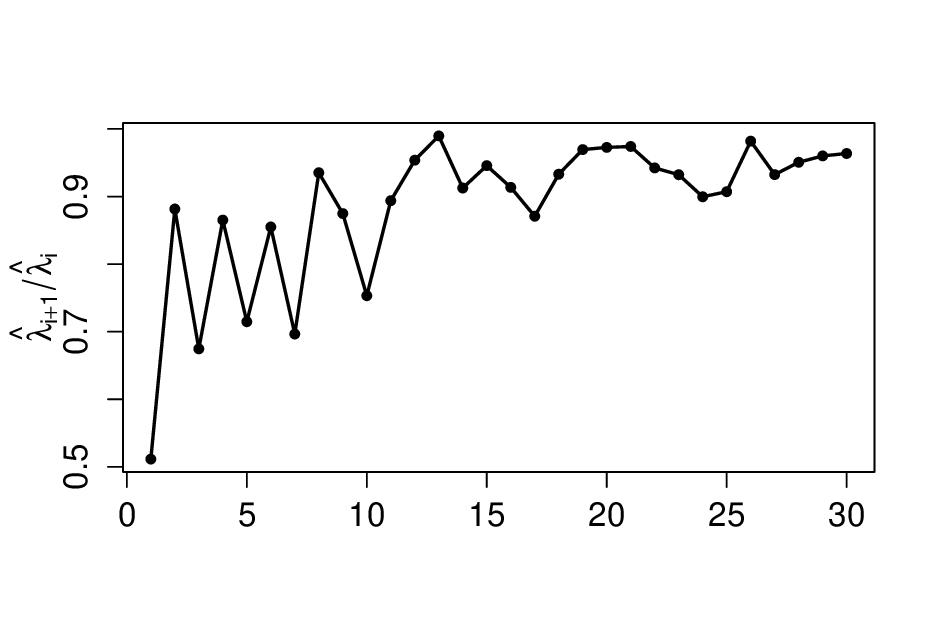}}
\subfigure[The estimated latent factor ]
{\includegraphics[width=0.8\linewidth]{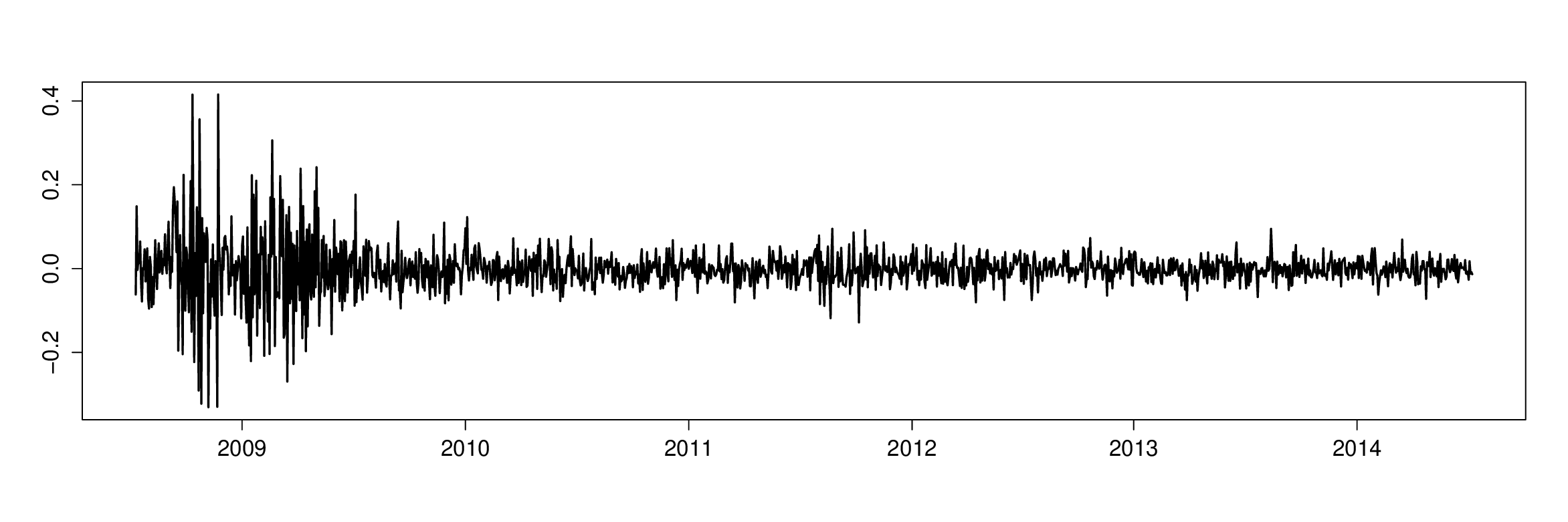}}     % 1500*500
\subfigure[S\&P $500$ returns]
{\includegraphics[width=0.8\linewidth]{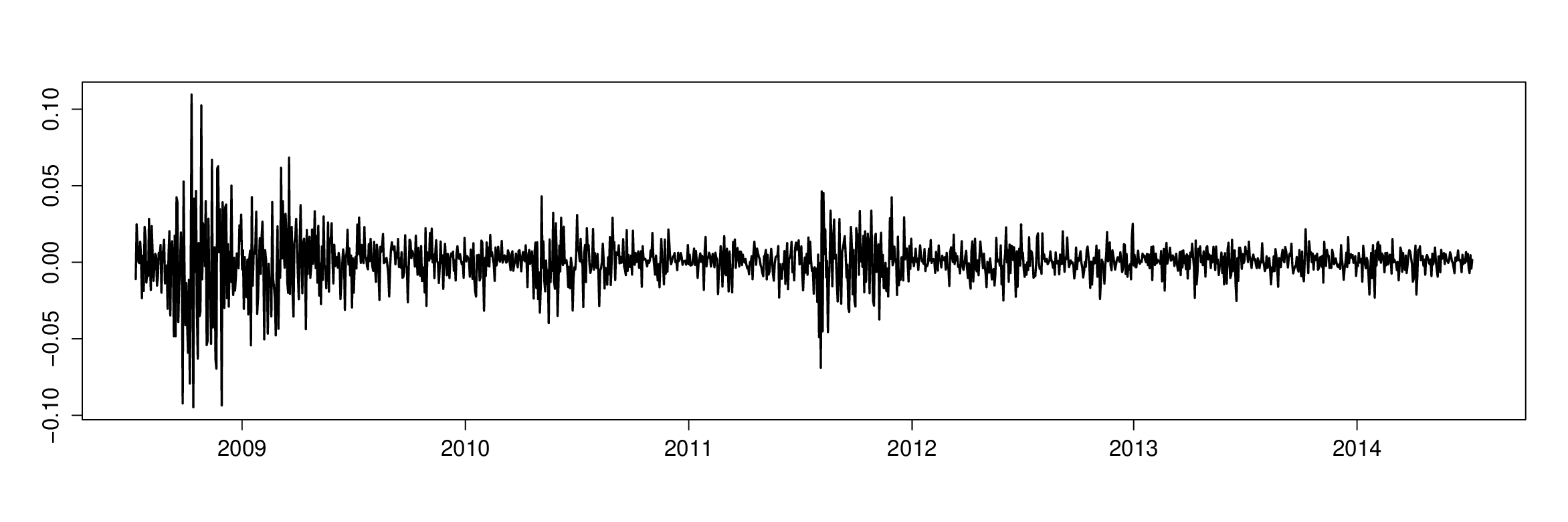}}
\caption{The estimated eigenvalues (multiplied by $10^6$), the ratio
of eigenvalues, the estimated latent factor and the S\&P$500$ returns in
14 July 2008 -- 11 July 2014.} \label{figure4}
\end{figure}

%\begin{figure}
%\centering \subfigure{\includegraphics[scale=0.36]{figure/acf1.eps}} %500*500
%\subfigure{\includegraphics[scale=0.36]{figure/acf2.eps}}
%\subfigure{\includegraphics[scale=0.36]{figure/acf3.eps}}
%\caption{The ACFs of the first three estimated factors.}
%\label{figure5}
%\end{figure}

\begin{figure}
\centering
\subfigure{\includegraphics[scale=0.5]{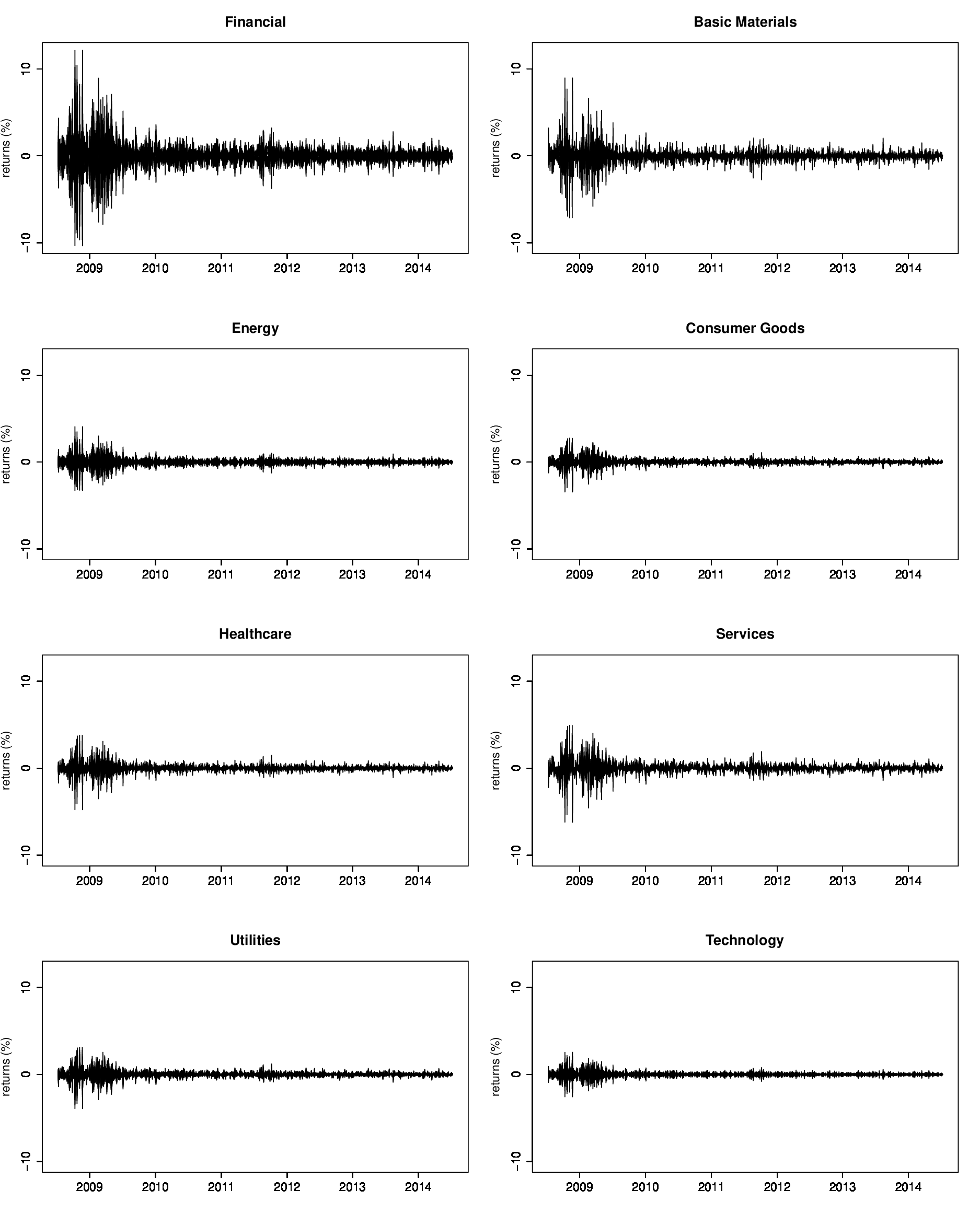}}
\caption{The estimated latent part ${\bf A}{\bf x}_t$ across
different sectors.} \label{figure6}
\end{figure}

To gain some appreciation of the latent factor, we divide the $123$
stocks into eight sectors: Financial, Basic Materials, Industrial
Goods, Consumer Goods, Healthcare, Services, Utilities and
Technology. We estimated the latent factor for each of those eight
sectors. Those estimated sector factors are plotted in
Fig.\ref{figure3}.  We observe that those estimated sector factors
behave differently for the different sectors. Especially the Basic
Materials sector exhibits the largest fluctuation. Consequently, we
may deduce that the oscillations, especially the sparks, of the
estimated factor in Fig.\ref{figure1}(c) are largely due to changes
in the Basic Materials sector. This is consistent with the relevant
economics and finance principles. Basic Materials sector includes
mainly the stocks of energy companies such as  oil, gas, coal et al.
The energy, especially oil, is the foundation for economic and
social development. Hence, the changes in oil price are often
considered as important events which underpin stock market
fluctuations, see, e.g. \cite{JonesKaul_1996} and
\cite{KilianPark_2009}. During January 2002 to December 2003,
international oil price had a huge increase. It rose $19\%$ from the
average in 2002. The 2003 invasion of Iraq marks a significant event
as Iraq possesses a significant portion of the global oil reserve.
Hence, the returns of the Basic Materials sector oscillate
dramatically during that period. Among other sectors, Industrial and
Consumer Goods have similar behaviors. However, the returns of both
the sectors have little changes around zero, thus they have little
contributions to the estimated factor. The same arguments hold for
the Utilities sector. Also note that the returns for the  Financial,
Healthcare, Services and Technology sectors are much less volatile
in comparison to that of the Basic Materials sector. We may conclude
that, the estimated factor mainly reflects the feature of stocks in
Basic Materials sector. The factor also contains some market
information about the Financial, Healthcare, Services and Technology
sectors, but less so on the Industrial Goods, Consumer Goods and
Utilities sectors.

We repeat the above exercise for another set of return data in
14 July 2008 -- 11 July 2014 from the 196 stocks contained in S\&P$500$.
Now $T=1510$ and $p=196$. The ratios of
$\widehat{\lambda}_{i+1}/\widehat{\lambda}_{i}$
shown in Fig.\ref{figure4}(b) indicate that there is still only one
latent factor, in addition to S\&P$500$. The estimated latent factor shown in
Fig.\ref{figure4}(c) fluctuated widely around 2009, which is
consistent with the pronounced decline the stock market due
to the global financial crisis.
While the latent factor process seems to resemble the returns
of S\&P$500$ (see  Figs.\ref{figure4}(c) and (d)), the two series are
orthogonal with each other (with the sample correlation coefficient equal to 0.00047).
The estimated factors for each of the eight sectors are
plotted in Fig.\ref{figure6}. In contrast to the findings in 2002-2008,
all the eight sectors contributed to the fluctuation around 2009, though
the financial sector was most predominant. The crisis caused by
the sharp down-turn of financial industry in early 2009 impacted all
sectors in the society.

\section*{Acknowledgements}
The authors would like to thank two referees for their helpful comments
and suggestions.

\section*{Appendix}

Throughout the Appendix, we use $C$s to denote generic uniformly
positive constants only depends on the parameters $C_i$s
 appear in the technical conditions which may be different in
different uses. Meanwhile, we denote ${\bf A}{\bf x}_t$ by
$\bzeta_t$. We first present the following lemmas which are used in
proofs of the propositions and theorems.

\begin{la}
\label{lemma1} Under Conditions \ref{as1}-\ref{as2}, $
\|T^{-1}\sum_{t=1}^T\{{\bf z}_t{\bf z}_t^\T-E({\bf z}_t{\bf
z}_t^\T)\}\|_F=O_p(mT^{-1/2}). $
\end{la}
\noindent\textbf{Proof}: For any $i,j=1,\ldots,m$, by Cauchy-Schwarz
inequality and Davydov inequality,
\begin{equation}
\begin{split}
&~E\bigg[\bigg|\frac{1}{T}\sum_{t=1}^T\{z_{i,t}z_{j,t}-E(z_{i,t}z_{j,t})\}\bigg|^2\bigg]\\
=&~\frac{1}{T^2}\sum_{t=1}^TE[\{z_{i,t}z_{j,t}-E(z_{i,t}z_{j,t})\}^2]\\
&+\frac{1}{T^2}\sum_{t_1\neq
t_2}E[\{z_{i,t_1}z_{j,t_1}-E(z_{i,t_1}z_{j,t_1})\}\{z_{i,t_2}z_{j,t_2}-E(z_{i,t_2}z_{j,t_2})\}]\\
\leq&~\frac{C}{T}+\frac{C}{T^2}\sum_{t_1\neq
t_2}\alpha(|t_1-t_2|)^{1-2/\gamma}\leq\frac{C}{T}+\frac{C}{T}\sum_{u=1}^{T}\alpha(u)^{1-2/\gamma}.\label{eq:a1}
\end{split}
\end{equation}
Then, $E\{\|T^{-1}\sum_{t=1}^T\{{\bf z}_t{\bf z}_t^\T-E({\bf
z}_t{\bf z}_t^\T)\}\|_F^2\}=O(m^2T^{-1})$ which implies the result.
$\hfill\Box$

\medskip

\noindent \textsc{Proof of Proposition \ref{pn1}}: Note that $
(\widehat{{\bf D}}-{\bf D})^\T=(T^{-1}\sum_{t=1}^T{\bf z}_t{\bf
z}_t^\T)^{-1}(T^{-1}\sum_{t=1}^T{\bf z}_t\beeta_t^\T) $ and
$\lambda_{\min}(T^{-1}\sum_{t=1}^T{\bf z}_t{\bf z}_t^\T)$ is bounded
away from zero with probability approaching one, which is implied by
Condition \ref{as3} and Lemma \ref{lemma1}, then $ \|\widehat{{\bf
D}}-{\bf D}\|_F=O_p(\|T^{-1}\sum_{t=1}^T{\bf z}_t\beeta_t^\T\|_F)$.
For each $i=1,\ldots,m$ and $j=1,\ldots,p$, from $\textrm{cov}({\bf
z}_t,\beeta_t)={\bf 0}$ and
%\[
%E\bigg(\bigg|\frac{1}{T}\sum_{t=1}^Tz_{i,t}\eta_{j,t}\bigg|^2\bigg)=
%E\bigg[\bigg|\frac{1}{T}\sum_{t=1}^T\{z_{i,t}\eta_{j,t}-E(z_{i,t}\eta_{j,t})\}\bigg|^2\bigg].
%\]
similar to (\ref{eq:a1}), we can obtain $
E\{(T^{-1}\sum_{t=1}^Tz_{i,t}\eta_{j,t})^2\} \leq CT^{-1}$. Then,
$E\{\|T^{-1}\sum_{t=1}^T{\bf z}_t\beeta_t^\T\|_F^2\}=O(pT^{-1})$.
Hence, $\|\widehat{{\bf D}}-{\bf D}\|_F=O_p(p^{1/2}T^{-1/2})$.
$\hfill \square$

\begin{la}
\label{lemma2} Under Conditions \ref{as1}-\ref{as2}, if $k=o(T)$,
then
\[
\|\widehat{\bSigma}_\eta(k)-\bSigma_\eta(k)\|_F=\|{\bf
D}-\widehat{{\bf D}}\|_F^2 J_{1,k}+\|{\bf D}-\widehat{{\bf D}}\|_F
J_{2,k}+J_{3,k}
\]
where $E(J_{1,k}^2)\leq Ckm^2(T-k)^{-1}+Cm^2\alpha(k)^{2-4/\gamma}$,
$E(J_{2,k}^2)\leq Ckpm(T-k)^{-1}+Cpm\alpha(k)^{2-4/\gamma}$ and
$E(J_{3,k}^2)\leq Ckp^2(T-k)^{-1}$.
\end{la}
\noindent\textbf{Proof}: For each $k=o(T)$,
\[
\begin{split}
\widehat{\bSigma}_\eta(k)-\bSigma_\eta(k)%&~\frac{1}{T-k}\sum_{t=1}^{T-k}\widehat{\beeta}_{t+k}\widehat{\beeta}_t^\T-E(\beeta_{t+k}\beeta_t^\T)\\
=&~\frac{1}{T-k}\sum_{t=1}^{T-k}(\widehat{\beeta}_{t+k}\widehat{\beeta}_t^\T-\beeta_{t+k}\beeta_t^\T)+\frac{1}{T-k}\sum_{t=1}^{T-k}\{\beeta_{t+k}\beeta_t^\T-E(\beeta_{t+k}\beeta_t^\T)\}\\
&~+\bar{\beeta}\bar{\beeta}^\T-\frac{1}{T-k}\sum_{t=1}^{T-k}\widehat{\beeta}_{t+k}\bar{\beeta}^\T-\frac{1}{T-k}\sum_{t=1}^{T-k}\bar{\beeta}\widehat{\beeta}_t^\T\\
=&~I_{1,k}+I_{2,k}+I_{3,k}+I_{4,k}+I_{5,k}.
\end{split}
\]
As
\[
\begin{split}
%&~\frac{1}{T-k}\sum_{t=1}^{T-k}(\widehat{\beeta}_{t+k}\widehat{\beeta}_t^\T-\beeta_{t+k}\beeta_t^\T)\\
%=&~\frac{1}{T-k}\sum_{t=1}^{T-k}(\widehat{\beeta}_{t+k}-\beeta_{t+k})(\widehat{\beeta}_t-\beeta_t)^\T
%+\frac{1}{T-k}\sum_{t=1}^{T-k}(\widehat{\beeta}_{t+k}-\beeta_{t+k})\beeta_t^\T\\
%&+\frac{1}{T-k}\sum_{t=1}^{T-k}\beeta_{t+k}(\widehat{\beeta}_t-\beeta_t)^\T\\
I_{1,k}=&~({\bf D}-\widehat{{\bf
D}})\bigg(\frac{1}{T-k}\sum_{t=1}^{T-k}{\bf z}_{t+k}{\bf
z}_t^\T\bigg)({\bf D}-\widehat{{\bf D}})^\T+({\bf D}-\widehat{{\bf
D}})\bigg(\frac{1}{T-k}\sum_{t=1}^{T-k}{\bf
z}_{t+k}\beeta_t^\T\bigg)\\
&+\bigg(\frac{1}{T-k}\sum_{t=1}^{T-k}\beeta_{t+k}{\bf
z}_t^\T\bigg)({\bf D}-\widehat{{\bf D}})^\T,\\
\end{split}
\]
then
\[
\begin{split}
%&~\bigg\|\frac{1}{T-k}\sum_{t=1}^{T-k}(\widehat{\beeta}_{t+k}\widehat{\beeta}_t^\T-\beeta_{t+k}\beeta_t^\T)\bigg\|_F\\
\|I_{1,k}\|_F\leq&~ \|{\bf D}-\widehat{{\bf
D}}\|_F^2\bigg\|\frac{1}{T-k}\sum_{t=1}^{T-k}{\bf z}_{t+k}{\bf
z}_t^\T\bigg\|_F+\|{\bf D}-\widehat{{\bf
D}}\|_F\bigg\|\frac{1}{T-k}\sum_{t=1}^{T-k}{\bf
z}_{t+k}\beeta_t^\T\bigg\|_F\\
&+\|{\bf D}-\widehat{{\bf
D}}\|_F\bigg\|\frac{1}{T-k}\sum_{t=1}^{T-k}\beeta_{t+k}{\bf
z}_t^\T\bigg\|_F.
\end{split}
\]

For any $i,j=1,\ldots,m$,
\[
\begin{split}
E\bigg\{\bigg(\frac{1}{T-k}\sum_{t=1}^{T-k}z_{i,t+k}z_{j,t}\bigg)^2\bigg\}\leq
2E\bigg(\bigg[\frac{1}{T-k}\sum_{t=1}^{T-k}\{z_{i,t+k}z_{j,t}-E(z_{i,t+k}z_{j,t})\}
\bigg]^2\bigg)+
2 \max_t \{E(z_{i,t+k}z_{j,t})\}^2.
\end{split}
\]
By Cauchy-Schwarz inequality and Davydov inequality,
\begin{equation}
\begin{split}
&~E\bigg(\bigg[\frac{1}{T-k}\sum_{t=1}^{T-k}\{z_{i,t+k}z_{j,t}-E(z_{i,t+k}z_{j,t})\}\bigg]^2\bigg)\\
=&~\frac{1}{(T-k)^2}\sum_{t=1}^{T-k}E[\{z_{i,t+k}z_{j,t}-E(z_{i,t+k}z_{j,t})\}^2]\\
&+\frac{1}{(T-k)^2}\sum_{t_1\neq
t_2}E[\{z_{i,t_1+k}z_{j,t_1}-E(z_{i,t_1+k}z_{j,t_1})\}\{z_{i,t_2+k}z_{j,t_2}-E(z_{i,t_2+k}z_{j,t_2})\}]\\
%\leq&~\frac{C}{T-k}+\frac{C}{(T-k)^2}\sum_{t_1<t_2}\alpha([t_2-(t_1+k)]^+)^{1-2/\gamma}\\
%\leq&~
%\frac{2C}{T-k}+\frac{16C}{(T-k)^2}\sum_{t_1<t_2}\alpha([t_2-(t_1+k)]^+)^{1-2/\gamma}\\
\leq&~\frac{C}{T-k}+\frac{Ck}{T-k}+\frac{Ck(k-1)}{(T-k)^2}+\frac{C}{T-k}\sum_{u=1}^{T-2k-1}\alpha(u)^{1-2/\gamma}.\label{eq:a2}\\
\end{split}
\end{equation}
and $\{E(z_{i,t+k}z_{j,t})\}^2\leq C\alpha(k)^{2-4/\gamma}$. Then, $
E[\{(T-k)^{-1}\sum_{t=1}^{T-k}z_{i,t+k}z_{j,t}\}^2]\leq
Ck(T-k)^{-1}+C\alpha(k)^{2-4/\gamma}$. Thus, $
E\{\|(T-k)^{-1}\sum_{t=1}^{T-k}{\bf z}_{t+k}{\bf z}_t^\T\|_F^2\}\leq
Ckm^2(T-k)^{-1}+Cm^2\alpha(k)^{2-4/\gamma}. $ By the same argument,
we can obtain $ E\{\|(T-k)^{-1}\sum_{t=1}^{T-k}{\bf
z}_{t+k}\beeta_t^\T\|_F^2\}\leq
Ckpm(T-k)^{-1}+Cpm\alpha(k)^{2-4/\gamma}$ and $
E\{\|(T-k)^{-1}\sum_{t=1}^{T-k}\beeta_{t+k}{\bf z}_t^\T\|_F^2\}\leq
Ckpm(T-k)^{-1}+Cpm\alpha(k)^{2-4/\gamma}$. Hence, $
\|I_{1,k}\|_F=\|{\bf D}-\widehat{{\bf D}}\|_F^2 J_{1,k}+\|{\bf
D}-\widehat{{\bf D}}\|_F J_{2,k} $ where $E(J_{1,k}^2)\leq
Ckm^2(T-k)^{-1}+Cm^2\alpha(k)^{2-4/\gamma}$ and $E(J_{2,k}^2)\leq
Ckpm(T-k)^{-1}+Cpm\alpha(k)^{2-4/\gamma}$. On the other hand,
similar to (\ref{eq:a2}), we can obtain $ E(\|I_{2,k}\|_F^2)\leq
Ckp^2(T-k)^{-1}$. For $I_{3,k}$, we have $E(\|I_{3,k}\|_F^2)\leq
E(\|\bar{\beeta}\|_2^4)$. By Jensen inequality and Davydov
inequality, $E(\|I_{3,k}\|_F^2)\leq Cp^2T^{-1}$. Following the same
way, we have both $E(\|I_{4,k}\|_F^2)$ and $E(\|I_{5,k}\|_F^2)$ can
be bounded by $Ckp^2(T-k)^{-1}$. Hence, we complete the proof.
$\hfill\Box$

%\begin{la}
%\label{lemma3} Under Condition \ref{as4},
%\[
%Cp^{(1-\delta)/2}\leq\|{\bf A}\|_{\min}\leq\|{\bf A}\|_2\leq
%C^{-1}p^{(1-\delta)/2}.
%\]
%\end{la}
%\noindent \textbf{Proof}: By the QR decomposition ${\bf A}={\bf
%Q}{\bf R}$, where ${\bf Q}$ is an orthogonal matrix such that ${\bf
%Q}^\T{\bf Q}=I_{r}$ and ${\bf R}$ is an upper triangular matrix.
%It is obvious that $\|{\bf
%A}\|_{\min}=\min_{i=1,\ldots,r}|R_{i,i}|$ and $\|{\bf
%A}\|_2=\max_{i=1,\ldots,r}|R_{i,i}|$ where $R_{i,i}$ is the
%$(i,i)$th component of ${\bf R}$. Let ${\bf R}=({\bf
%r}_1,\ldots,{\bf r}_{r})$, then $\|{\bf a}_i-\sum_{j\neq
%i}\theta_j{\bf a}_j\|_2^2=\|{\bf r}_i-\sum_{j\neq i}\theta_j{\bf
%r}_j\|_2^2. $ For $i=1$, setting $\theta_j=0$ for $j\neq 1$ and
%noting Condition \ref{as4}, we have $|R_{1,1}|^2=\|{\bf
%r}_1\|_2^2=\|{\bf r}_1-\sum_{j\neq1}\theta_j{\bf r}_j\|_2^2\asymp
%p^{1-\delta}$.
% Let
%$\theta_1=R_{1,2}R_{1,1}^{-1}$, $\theta_{2}=1$ and $\theta_{j}=0$
%for each $j\geq 3$, then $|R_{2,2}|^2=\|{\bf
%r}_2-\sum_{j\neq2}\theta_j{\bf r}_j\|_2^2\asymp p^{1-\delta}$.
%Repeat the same argument, we have $|R_{i,i}|^2\asymp p^{1-\delta}$
%for any $i=1,\ldots,r$. Hence, we complete the proof of this
%lemma. $\hfill\Box$

\begin{la}
\label{lemma4} Under Condition \ref{as5}, for $k=1,\ldots,\bar{k}$,
\[
\|\bSigma_\eta(k)\|_2\leq Cp^{1-\delta}+C\kappa_2.
\]
\end{la}
\noindent  \textbf{Proof}: Note that $ \bSigma_\eta(k)={\bf
A}\bSigma_x(k){\bf A}^\T+\bSigma_{x\varepsilon}(k)$, then $
\|\bSigma_\eta(k)\|_2\leq \|{\bf
A}\|_2^2\|\bSigma_x(k)\|_2+\|\bSigma_{x\varepsilon}(k)\|_2. $ From
 Condition \ref{as5}, we complete the proof.
$\hfill\Box$

\begin{la}
\label{lemma5} Under Conditions \ref{as1}-\ref{as5},
\[
\|\widehat{{\bf M}}-{\bf
M}\|_2=O_p\{(p^{1-\delta}+\kappa_2)pT^{-1/2}+p^2T^{-1}\}.
\]
\end{la}

\noindent  \textbf{Proof}: Note that
\[
\begin{split}
\|\widehat{{\bf M}}-{\bf M}\|_2%\leq&~
%\sum_{k=1}^{\bar{k}}\|\widehat{\Sigma}_\eta(k)\widehat{\Sigma}_\eta(k)^\T-\Sigma_\eta(k)\Sigma_\eta(k)^\T\|_2\\
\leq&~\sum_{k=1}^{\bar{k}}\|\widehat{\bSigma}_\eta(k)-\bSigma_\eta(k)\|_2^2+2\sum_{k=1}^{\bar{k}}\|\bSigma_\eta(k)\|_2\|\widehat{\bSigma}_\eta(k)-\bSigma_\eta(k)\|_2=I_1+I_2.\\
\end{split}
\]
By Lemmas \ref{lemma2} and \ref{lemma4}, we can obtain
\[
\begin{split}
I_1\leq&~ 3\|{\bf D}-\widehat{{\bf
D}}\|_F^4\sum_{k=1}^{\bar{k}}J_{1,k}^2+3\|{\bf D}-\widehat{{\bf
D}}\|_F^2\sum_{k=1}^{\bar{k}}J_{2,k}^2+3\sum_{k=1}^{\bar{k}}J_{3,k}^2=O_p(p^2T^{-1})
\end{split}
\]
and
\[
\begin{split}
I_2\leq&~2\bigg\{\|{\bf D}-\widehat{{\bf
D}}\|_F^2\sum_{k=1}^{\bar{k}}J_{1,k}+\|{\bf D}-\widehat{{\bf
D}}\|_F\sum_{k=1}^{\bar{k}}J_{2,k}+\sum_{k=1}^{\bar{k}}J_{3,k}\bigg\}\sup_{1\leq
k\leq\bar{k}}\|\bSigma_\eta(k)\|_2\\
=&~O_p\{(p^{1-\delta}+\kappa_2)pT^{-1/2}\}.
\end{split}
\]
Hence, we complete the proof. $\hfill\Box$

\begin{la}
\label{lemma6}Under Condition \ref{as5},
\[
\lambda_{r}({\bf M})\geq \left\{ \begin{aligned}
         Cp^{2(1-\delta)},&~~\textrm{if}~\kappa_2=o(p^{1-\delta});\\
         C\kappa_1^2,&~~\textrm{if}~
p^{1-\delta}=o(\kappa_1).
        \end{aligned} \right.
\]
\end{la}

\noindent \textbf{Proof}: From (\ref{eq:m}), we know
\[
\lambda_{r}({\bf
M})=\lambda_{\min}\bigg[\sum_{k=1}^{\bar{k}}\{\bSigma_x(k){\bf
A}^\T+\bSigma_{x\varepsilon}(k)\}\{\bSigma_x(k){\bf
A}^\T+\bSigma_{x\varepsilon}(k)\}^\T\bigg].
\]
For each $k=1,\ldots,\bar{k}$,
\[
\begin{split}
&~\lambda_{\min}[\{\bSigma_x(k){\bf
A}^\T+\bSigma_{x\varepsilon}(k)\}\{\bSigma_x(k){\bf
A}^\T+\bSigma_{x\varepsilon}(k)\}^\T]\\
 \asymp&~ \left\{
\begin{aligned}
         \lambda_{\min}\{\bSigma_x(k)\bSigma_x(k)^\T\},&~~\textrm{if}~\lambda_{\max}\{\bSigma_{x\varepsilon}(k)\bSigma_{x\varepsilon}(k)^\T\}=o(\lambda_{\min}\{\bSigma_x(k)\bSigma_x(k)^\T\}); \\
         \lambda_{\min}\{\bSigma_{x\varepsilon}(k)\bSigma_{x\varepsilon}(k)^\T\},&~~\textrm{if}~\lambda_{\max}\{\bSigma_{x}(k)\bSigma_{x}(k)^\T\}=o(\lambda_{\min}\{\bSigma_{x\varepsilon}(k)\bSigma_{x\varepsilon}(k)^\T\}).\\
        \end{aligned} \right.\\
\end{split}
\]
Notice Condition \ref{as5}, then
\[
\lambda_{\min}[\{\bSigma_x(k){\bf
A}^\T+\bSigma_{x\varepsilon}(k)\}\{\bSigma_x(k){\bf
A}^\T+\bSigma_{x\varepsilon}(k)\}^\T]\geq \left\{
\begin{aligned}
         Cp^{2(1-\delta)},&~~\textrm{if}~\kappa_2=o(p^{1-\delta}); \\
         C\kappa_1^2,&~~\textrm{if}~p^{1-\delta}=o(\kappa_1).\\
        \end{aligned} \right.\\
\]
Hence, we complete the proof.$\hfill\Box$

\medskip

\noindent \textsc{Proof of Theorem \ref{tm1}}: %From Lemma
%\ref{lemma5},
%\[
%\begin{split}
%\|\widehat{{\bf M}}-{\bf M}\|_2=&~\left\{
%     \begin{aligned}
%       O_p(p^{2-\delta}T^{-1/2}+p^2T^{-1}),  &~~\hbox{if $\kappa_2p^{-(1-\delta)/2}=o(1)$;} \\
%       O_p(p^{(3-\delta)/2}\kappa_2T^{-1/2}+p^2T^{-1}), &~~\hbox{if $\kappa_1^{-1}p^{(1-\delta)/2}=o(1)$.}
%     \end{aligned}
%   \right.
%\end{split}
%\]
By Lemma \ref{lemma6}, $\|\widehat{{\bf M}}-{\bf
M}\|_2=o_p\{\lambda_{r}({\bf M})\} $ provided that either case (i)
$\kappa_2=o(p^{1-\delta})$ and $p^{2\delta}T^{-1}=o(1)$ or (ii)
$p^{1-\delta}=o(\kappa_1)$ and $\kappa_1^{-2}\kappa_2pT^{-1/2}=o(1)$
hold. By Lemma 3 of \cite{LamYaoBathia_Biometrika_2011}, and using
the same argument of the proof of Theorem 1 in their paper, %we have
\[
\begin{split}
\|\widehat{{\bf A}}-{\bf A}\|_2 =&~\left\{
     \begin{aligned}
       O_p(p^\delta
T^{-1/2}),&~~\hbox{if $\kappa_2=o(p^{1-\delta})$ and
$p^{2\delta} T^{-1}=o(1)$;} \\
       O_p(\kappa_1^{-2}\kappa_2pT^{-1/2}),&~~\hbox{if $p^{1-\delta}=o(\kappa_1)$ and
$\kappa_1^{-2}\kappa_2pT^{-1/2}=o(1)$}.
     \end{aligned}
   \right.
\end{split}
\]
Hence, we complete the proof. $\hfill\Box$

\medskip

\noindent \textsc{Proof of Theorem \ref{tm2}}: Note that
\[
\begin{split}
\widehat{{\bf A}}\widehat{{\bf x}}_t-{\bf A}{\bf
x}_t=&~\widehat{{\bf A}}\wh\bA^\T{\bf A}{\bf x}_t+\widehat{{\bf
A}}\wh\bA^\T\bvarepsilon_t-{\bf A}{\bf x}_t+\widehat{{\bf A}}\wh\bA^\T(\widehat{\beeta}_t-\beeta_t)\\
=&~(\widehat{{\bf A}}\wh\bA^\T-{\bf A}\bA^\T){\bf A}{\bf
x}_t+\widehat{{\bf A}}(\widehat{{\bf A}}-{\bf
A})^\T\bvarepsilon_t+\widehat{{\bf
A}}\bA^\T\bvarepsilon_t+\widehat{{\bf A}}\wh\bA^\T(\widehat{\beeta}_t-\beeta_t)\\
=&~I_1+I_2+I_3+I_4.
\end{split}
\]
For $I_1$, $\|I_1\|_2\leq 2\|\widehat{{\bf A}}-{\bf A}\|_2\|{\bf
A}{\bf x}_t\|_2\leq O_p(p^{1/2}\|\widehat{{\bf A}}-{\bf A}\|_2)$.
For $I_2$, $ \|I_2\|_2\leq \|\widehat{{\bf A}}-{\bf A}\|_2\|
\bvarepsilon_t\|_2=O_p(p^{1/2}\|\widehat{{\bf A}}-{\bf A}\|_2). $
For $I_3$,
 as $
E(\|I_3\|_2^2)=\sum_{i=1}^{r}E\{({\bf a}_i^\T\bvarepsilon_t)^2\}\leq
r\lambda_{\max}(\bSigma_{\varepsilon}) $, then $ I_3=O_p(1). $ For
$I_4$, by Proposition \ref{pn1}, $ \|I_4\|_2\leq \|\widehat{{\bf
D}}-{\bf D}\|_2\|{\bf z}_t\|_2=O_p(p^{1/2}T^{-1/2}). $ Hence, $
p^{-1/2}\|\widehat{{\bf A}}\widehat{{\bf x}}_t-{\bf A}{\bf
x}_t\|_2\leq O_p(\|\widehat{{\bf A}}-{\bf A}\|_2+p^{-1/2}+T^{-1/2}).
$ $\hfill\Box$

\medskip

\noindent \textsc{Proof of Theorem \ref{tm3}}: Let
$\bSigma_{\zeta}(k)=(T-k)^{-1}\sum_{t=1}^{T-k}\textrm{cov}(\bzeta_{t+k},\bzeta_t)$,
then $\bSigma_{\zeta}(k)={\bf A}\bSigma_x(k){\bf A}^\T$. Note that
\[
\begin{split}
\textrm{tr}\{\bSigma_{\zeta}(1)^\T(I_p-\widehat{{\bf A}}\wh {\bf
A}^\T)\bSigma_{\zeta}(1)\}=&~\textrm{tr}\{\bSigma_x(1)^\T(I_{r}-{\bf
A}^\T\widehat{{\bf A}}\wh{\bf A}^\T{\bf A})\bSigma_x(1)\}\\
\geq&~\textrm{tr}(I_{r}-{\bf A}^\T\widehat{{\bf A}}\widehat{{\bf
A}}^\T{\bf A})\lambda_{\min}\{\bSigma_x(1)\bSigma_x(1)^\T\}\\
=&~r\{D(\mathcal{M}(\widehat{{\bf A}}),\mathcal{M}({\bf
A}))\}^2\lambda_{\min}\{\bSigma_x(1)\bSigma_x(1)^\T\}.
\end{split}
\]
By Condition \ref{as5},
\[
\textrm{tr}\{\bSigma_{\zeta}(1)^\T(I_p-\widehat{{\bf
A}}\widehat{{\bf A}}^\T)\bSigma_{\zeta}(1)\}\geq
Cp^{2(1-\delta)}\{D(\mathcal{M}(\widehat{{\bf A}}),\mathcal{M}({\bf
A}))\}^2.
\]
At the same time,
\[
\begin{split}
&~\textrm{tr}\{\bSigma_{\zeta}(1)^\T(I_p-\widehat{{\bf
A}}\widehat{{\bf
A}}^\T)\bSigma_{\zeta}(1)\}-\textrm{tr}\{\bSigma_{\zeta}(1)^\T(I_p-{\bf
A}{\bf A}^\T)\bSigma_{\zeta}(1)\}\\
=&~\textrm{tr}\{{\bf A}\bSigma_{x}(1)^\T{\bf A}^\T({\bf A}{\bf
A}^\T-\widehat{{\bf
A}}\widehat{{\bf A}}^\T){\bf A}\bSigma_{x}(1){\bf A}^\T\}\\
\leq&~\lambda_{\max}\{{\bf A}^\T({\bf A}{\bf A}^\T-\widehat{{\bf
A}}\widehat{{\bf A}}^\T){\bf
A}\}\textrm{tr}\{\bSigma_{x}(1)\bSigma_{x}(1)^\T\}\\
\leq&~Cp^{2(1-\delta)}\|{\bf A}^\T({\bf A}{\bf A}^\T-\widehat{{\bf
A}}\widehat{{\bf A}}^\T){\bf A}\|_2.
\end{split}
\]
Note that $\textrm{tr}\{\bSigma_{\zeta}(1)^\T(I_p-{\bf A}{\bf
A}^\T)\bSigma_{\zeta}(1)\}=0$, then
\[
\{D(\mathcal{M}(\widehat{{\bf A}}),\mathcal{M}({\bf A}))\}^2\leq
C\|{\bf A}^\T({\bf A}{\bf A}^\T-\widehat{{\bf A}}\widehat{{\bf
A}}^\T){\bf A}\|_2.
\]

On the other hand, we have the following two inequality,
\[
\begin{split}
\textrm{tr}\{\bSigma_{\zeta}(1)^\T(I_p-\widehat{{\bf
A}}\widehat{{\bf A}}^\T)\bSigma_{\zeta}(1)\}\leq&~
r\{D(\mathcal{M}(\widehat{{\bf A}}),\mathcal{M}({\bf
A}))\}^2\lambda_{\max}\{\bSigma_{x}(1)\bSigma_x(1)^\T\}\\
\leq&~Cp^{2(1-\delta)}\{D(\mathcal{M}(\widehat{{\bf
A}}),\mathcal{M}({\bf A}))\}^2
\end{split}
\]
and
\[
\begin{split}
\textrm{tr}\{\bSigma_{\zeta}(1)^\T(I_p-\widehat{{\bf
A}}\widehat{{\bf A}}^\T)\bSigma_{\zeta}(1)\}\geq&~
\lambda_{\min}\{\bSigma_x(1)\bSigma_x(1)^\T\}\textrm{tr}\{{\bf
A}^\T({\bf A}{\bf A}^\T-\widehat{{\bf A}}\widehat{{\bf A}}^\T){\bf
A}\}\\
\geq&~Cp^{2(1-\delta)}\|{\bf A}^\T({\bf A}{\bf A}^\T-\widehat{{\bf
A}}\widehat{{\bf A}}^\T){\bf
A}\|_2.%\\
%\geq&~\lambda_{\min}\{\Sigma_f(1)\Sigma_f(1)^\T\}\{\|\widehat{{\bf
%Q}}-{\bf Q}\|_2^2-\|{\bf Q}^\T(\widehat{{\bf Q}}-{\bf Q})\|_2^2\}.
\end{split}
\]
Hence,
\[
\{D(\mathcal{M}(\widehat{{\bf A}}),\mathcal{M}({\bf
A}))\}^2\asymp\|{\bf A}^\T({\bf A}{\bf A}^\T-\widehat{{\bf
A}}\widehat{{\bf A}}^\T){\bf A}\|_2.
\]
Note that \[ {\bf A}^\T({\bf A}{\bf A}^\T-\widehat{{\bf
A}}\widehat{{\bf A}}^\T){\bf A}=-{\bf A}^\T({\bf A}-\widehat{{\bf
A}})({\bf A}-\widehat{{\bf A}})^\T{\bf A}+({\bf A}-\widehat{{\bf
A}})^\T({\bf A}-\widehat{{\bf A}}), \] then we complete the proof.
$\hfill\Box$

\medskip

\noindent \textsc{Proof of Theorem \ref{tm4}}: As
$(p^{1-\delta}+\kappa_2)pT^{-1/2}\log T\rightarrow0$, then
$\|\widehat{{\bf M}}-{\bf M}\|_2=o_p\{\lambda_{r}({\bf M})\}$. Then
$\sup_{j=1,\ldots,p}|\widehat{\lambda}_j-\lambda_j({\bf
M})|\leq\|\widehat{{\bf M}}-{\bf M}\|_2=o_p\{\lambda_{r}({\bf
M})\}$. For any $j<r$,
\[
\frac{\widehat{\lambda}_{j+1}+(p^{1-\delta}+\kappa_2)pT^{-1/2}\log
T}{\widehat{\lambda}_j+(p^{1-\delta}+\kappa_2)pT^{-1/2}\log
T}\xrightarrow{p}C>0.
\]
For any $j>r$, note that $\|\widehat{{\bf M}}-{\bf
M}\|_2=o_p\{(p^{1-\delta}+\kappa_2)pT^{-1/2}\log T\}$ which implies
that
$|\widehat{\lambda}_j|=o_p\{(p^{1-\delta}+\kappa_2)pT^{-1/2}\log
T\}$, then
\[
\frac{\widehat{\lambda}_{j+1}+(p^{1-\delta}+\kappa_2)pT^{-1/2}\log
T}{\widehat{\lambda}_j+(p^{1-\delta}+\kappa_2)pT^{-1/2}\log
T}\xrightarrow{p}1>0.
\]
On the other hand,
\[
\frac{\widehat{\lambda}_{r+1}+(p^{1-\delta}+\kappa_2)pT^{-1/2}\log
T}{\widehat{\lambda}_r+(p^{1-\delta}+\kappa_2)pT^{-1/2}\log
T}\xrightarrow{p}0.
\]
Hence, the criterion implies a consistent estimator of $r$.
$\hfill\Box$

\medskip

\noindent \textsc{Proof of Proposition \ref{pn2}}: Following the
proof of Lemma \ref{lemma1}, $\|T^{-1}\sum_{t=1}^T\{{\bf z}_t{\bf
w}_t^\T-E({\bf z}_t{\bf w}_t^\T)\}\|_F=O_p(m^{1/2}q^{1/2}T^{-1/2})$.
Note that $\textrm{rank}({\bf R})=m$ and Condition \ref{as7}, it
yields $\lambda_{\min}(T^{-1}\sum_{t=1}^T{\bf z}_t{\bf w}_t^\T{\bf
R}^\T)$ is bounded away from zero with probability approaching one.
Hence, following the proof of Proposition \ref{pn1}, we can obtain
the result. $\hfill\Box$

\bigskip

\noindent  \textsc{Proof of Proposition \ref{pnn5}}: For each
$i=1,\ldots,p$,
\[
\widehat{{\bf d}}_i-{\bf d}_i=\bigg(\frac{1}{T}\sum_{t=1}^T{\bf
z}_t{\bf
z}_t^\T\bigg)^{-1}\bigg(\frac{1}{T}\sum_{t=1}^T\eta_{i,t}{\bf
z}_t\bigg)+\bigg(\frac{1}{T}\sum_{t=1}^T{\bf z}_t{\bf
z}_t^\T\bigg)^{-1}\bigg(\frac{1}{T}\sum_{t=1}^Te_{i,t}{\bf
z}_t\bigg).
\]
Then,
\[
\|\widehat{{\bf d}}_i-{\bf
d}_i\|_2\lambda_{\min}\bigg(\frac{1}{T}\sum_{t=1}^T{\bf z}_t{\bf
z}_t^\T\bigg)\leq \bigg\|\frac{1}{T}\sum_{t=1}^T\eta_{i,t}{\bf
z}_t\bigg\|_2+\bigg\|\frac{1}{T}\sum_{t=1}^Te_{i,t}{\bf
z}_t\bigg\|_2.
\]
Note that $E(\bzeta_t|{\bf u}_t)={\bf 0}$ and $E(\bvarepsilon_t|{\bf
u}_t)={\bf 0}$, we have $ \|{T}^{-1}\sum_{t=1}^T\eta_{i,t}{\bf
z}_t\|_2=O_p(m^{1/2}T^{-1/2}) $ and $
\|{T}^{-1}\sum_{t=1}^Te_{i,t}{\bf
z}_t\|_2=\|{T}^{-1}\sum_{t=1}^TE(e_{i,t}{\bf
z}_t)\|_2+O_p(m^{1/2}T^{-1/2}), $ where $O_p(m^{1/2}T^{-1/2})$s are
uniformly for $i=1,\ldots,p$. On the other hand, $ \|E(e_{i,t}{\bf
z}_t)\|_2^2=O(m^{1-2\lambda}). $ Thus, we complete the
proof.\hfill$\Box$

\medskip

\noindent  \textsc{Proof of Theorem \ref{tm5}}: Let ${\bf
z}=(l_1({\bf u}),\ldots,l_m({\bf u}))^\T$. For each $i=1,\ldots,p$,
\[
\begin{split}
\widehat{g}_i({\bf u})-g_i({\bf u})%&~{\bf
%z}^\T\bigg(\frac{1}{T}\sum_{t=1}^T{\bf z}_t{\bf
%z}_t^\T\bigg)^{-1}\bigg(\frac{1}{T}\sum_{t=1}^Ty_{i,t}{\bf
%z}_t\bigg)-g_i({\bf u})\\
=&~{\bf z}^\T\bigg(\frac{1}{T}\sum_{t=1}^T{\bf z}_t{\bf
z}_t^\T\bigg)^{-1}\bigg\{\frac{1}{T}\sum_{t=1}^T{\bf
z}_t(e_{i,t}+\zeta_{i,t}+\varepsilon_{i,t})\bigg\}-e_i
\end{split}
\]
where $g_i({\bf u})={\bf d}_i^{\T}{\bf z}+e_i$. Hence,
\[
\begin{split}
\int_{{\bf u}\in\mathcal {U}}|\widehat{g}_i({\bf u})-g_i({\bf
u})|^2~d{\bf u}\leq&~2\bigg\{\frac{1}{T}\sum_{t=1}^T{\bf
z}_t^\T(e_{i,t}+\zeta_{i,t}+\varepsilon_{i,t})\bigg\}\bigg(\frac{1}{T}\sum_{t=1}^T{\bf
z}_t{\bf z}_t^\T\bigg)^{-1}\\
&~~~~~~\times\bigg(\int_{{\bf u}\in\mathcal{U}}{\bf z}{\bf
z}^\T~d{\bf u}\bigg)\bigg(\frac{1}{T}\sum_{t=1}^T{\bf z}_t{\bf
z}_t^\T\bigg)^{-1}\bigg\{\frac{1}{T}\sum_{t=1}^T{\bf
z}_t(e_{i,t}+\zeta_{i,t}+\varepsilon_{i,t})\bigg\}\\
&+Cm^{-2\lambda}.
\end{split}
\]
Let $p({\bf u}_t)$ be the density function of ${\bf u}_t$ and pick
${\bf v}^\T$ such that $ \lambda_{\max}(\int_{{\bf
u}\in\mathcal{U}}{\bf z}{\bf z}^\T~d{\bf u})=\int_{{\bf
u}\in\mathcal{U}}{\bf v}^\T{\bf z}{\bf z}^\T{\bf v}~d{\bf u}, $ by
Condition \ref{as12},
\[
{\bf v}^\T E({\bf z}_t{\bf z}_t^\T){\bf v}=\int_{{\bf
u}_t\in\mathcal{U}}{\bf v}^\T{\bf z}_t{\bf z}_t^\T{\bf v} p({\bf
u}_t)~d{\bf u}_t\geq C \int_{{\bf u}_t\in\mathcal{U}}{\bf v}^\T{\bf
z}_t{\bf z}_t^\T{\bf v} ~d{\bf u}_t=C\lambda_{\max}\bigg(\int_{{\bf
u}\in\mathcal{U}}{\bf z}{\bf z}^\T~d{\bf u}\bigg).
\]
From Condition \ref{as9}, we know $ \lambda_{\max}(\int_{{\bf
u}\in\mathcal{U}}{\bf z}{\bf z}^\T~d{\bf u})\leq C $ which implies
\[
\int_{{\bf u}\in\mathcal {U}}|\widehat{g}_i({\bf u})-g_i({\bf
u})|^2~d{\bf u}\leq O_p(mT^{-1})+O(m^{-2\lambda}).
\]
The terms $O_p(mT^{-1})$ and $O(m^{-2\lambda})$ are uniformly for
$i=1,\ldots,p$, thus we complete the proof. \hfill$\Box$

\begin{la}
\label{lemma9} For nonlinear regression model (\ref{eq:nonlinear}),
under Conditions \ref{as1}-\ref{as2}, \ref{as13} and \ref{as11}, if
$k=o(T)$, then
\[
\|\widehat{\bSigma}_\eta(k)-\bSigma_\eta(k)\|_F=\|{\bf
D}-\widehat{{\bf D}}\|_F^2 J_{1,k}+\|{\bf D}-\widehat{{\bf D}}\|_F
J_{2,k}+J_{3,k}
\]
where $E(J_{1,k}^2)\leq Ckm^2(T-k)^{-1}+Cm^2\alpha(k)^{2-4/\gamma}$,
$E(J_{2,k}^2)\leq Ckpm(T-k)^{-1}+Cpm\alpha(k)^{2-4/\gamma}$ and
$E(J_{3,k}^2)\leq
Ckp^2(T-k)^{-1}+Cp^2m^{-2\lambda}\alpha(k)^{2-4/\gamma}$.
\end{la}

\noindent\textbf{Proof}:
Noting $\sup_t\|{\bf e}_t\|_\infty=O(m^{-\lambda})$, similar to
Lemma \ref{lemma2}, we can obtain the result. $\hfill\Box$

\begin{la}
\label{lemma10} Under Conditions \ref{as1}-\ref{as2}, \ref{as5},
\ref{as13}-\ref{as11}, if $mT^{-1/2}=o(1)$, $\bar{k}T^{-1/2}=o(1)$
and $\lambda\geq1$, then
\[
\|\widehat{{\bf M}}-{\bf
M}\|_2=O_p\{(p^{1-\delta}+\kappa_2)p[(\bar{k}^{3/2}+m)T^{-1/2}+m^{1-\lambda}]\}+O_p\{p^2[(\bar{k}^2+m^2)T^{-1}+m^{2-2\lambda}]\}.
\]
\end{la}
\noindent \textbf{Proof}: Note that $ \|\widehat{{\bf M}}-{\bf
M}\|_2\leq
\sum_{k=1}^{\bar{k}}\{\|\widehat{\bSigma}_{\eta}(k)-\bSigma_{\eta}(k)\|_2^2+2\|\bSigma_\eta(k)\|_2\|\widehat{\bSigma}_\eta(k)-\bSigma_\eta(k)\|_2\}.
$ By Lemma \ref{lemma9}, we complete the proof. $\hfill\Box$

\medskip

\noindent  \textsc{Proof of Theorem \ref{tmn5}}: Note that
$m=O(T^{1/(2\lambda+1)})$, then
\[
\|\widehat{{\bf M}}-{\bf
M}\|_2=O_p\{(p^{1-\delta}+\kappa_2)p(\bar{k}^{3/2}T^{-1/2}+T^{(1-\lambda)(2\lambda+1)})+p^2(\bar{k}^2T^{-1}+T^{(2-2\lambda)/(2\lambda+1)})\}.
\]
Similar to the proof of Lemma \ref{lemma6}, we have
\[
\lambda_{r}({\bf M})\geq\left\{
                             \begin{aligned}
                               C\bar{k}p^{2(1-\delta)}, &~~ \hbox{if $\kappa_2=o(p^{1-\delta})$;} \\
                               C\bar{k}\kappa_1^2, &~~ \hbox{if $p^{1-\delta}=o(\kappa_1)$.}
                             \end{aligned}
                           \right.
\]
Then, by Lemma \ref{lemma10}, $\|\widehat{{\bf M}}-{\bf
M}\|_2=o_p\{\lambda_{r}({\bf M})\} $ provided that either (i)
$\kappa_2=o(p^{1-\delta})$ and
$p^{2\delta}[\bar{k}T^{-1}+T^{(2-2\lambda)/(2\lambda+1)}]=o(1)$ or
(ii) $p^{1-\delta}=o(\kappa_1)$ and
$p^{2}\kappa_2^2\kappa_1^{-4}[\bar{k}T^{-1}+T^{(2-2\lambda)/(2\lambda+1)}]=o(1)$
hold. Using the same argument of the proof of Theorem \ref{tm1}, we
obtain the result. $\hfill\Box$

\medskip

\noindent  \textsc{Proof of Theorem \ref{tmn6}}: Following the
arguments of the proof of Theorem \ref{tm2}, we can construct the
result. $\hfill\Box$

\newpage

\newpage

\end{document}